\documentclass[hidelinks,12pt]{article}

\usepackage{natbib}

\usepackage{graphicx}
\usepackage{subfigure}
\usepackage{amsmath,amsthm,amsopn,amstext,amscd,amsfonts,amssymb}
\usepackage{bbm}
\usepackage{bm}
\usepackage{multirow} 
\usepackage{setspace}
\usepackage{hyperref}
\usepackage{mathrsfs}
\usepackage{mdwlist}
\usepackage{rotating}
\usepackage{arydshln}
\usepackage{pdflscape}
\usepackage{enumitem}
\usepackage{color}

\hypersetup{
    colorlinks=false,
    pdfborder={0 0 0},
}

\makeindex

\newcommand{\field}[1]{\mathbb{#1}}
\newcommand{\R}{\field{R}}   
\newcommand{\N}{\field{N}}
\newcommand{\E}{\field{E}}

\newcommand{\ts}{^{\sf T}} 
\newcommand{\In}{\bm{\mathcal{I}}}
\newcommand{\eb}{ \bm\epsilon }
\newcommand*{\defeq}{\stackrel{\text{def}}{=}}
\newcommand{\Sl}{\textbf{S}_{\hat{\bm\lambda}}}
\newcommand{\I}{\textbf{I}}
\newcommand{\muz}{\bm\mu_{\textbf{z}}}
\newcommand{\Al}{\textbf{A}_{\bm\lambda}}
\newcommand{\beq}{\begin{equation}}
\newcommand{\eeq}{\end{equation}}

\def\bbeta{\boldsymbol\beta}
\def\balpha{\boldsymbol\alpha}
\def\bdelta{\boldsymbol\delta}
\theoremstyle{Conjecture}
\theoremstyle{example}
\theoremstyle{remark}
\theoremstyle{lemma}
\theoremstyle{definition}
\theoremstyle{corol}
\theoremstyle{proposition}
\theoremstyle{condition}

\def\bm{\boldsymbol}

\begin{document}


\title{Copula based generalized additive models with non-random sample selection}
\author{
Ma{\l}gorzata Wojty\'s \vspace{-1ex}\\ \small Centre for Mathematical Sciences \vspace{-1ex}\\ \small Plymouth University  \vspace{-1ex}\\ \small Drake Circus, Plymouth PL4 8AA, U.K.
\and
Giampiero Marra  \vspace{-1ex}\\ \small Department of Statistical Science  \vspace{-1ex}\\ \small University College London  \vspace{-1ex}\\ \small Gower Street, London WC1E 6BT, U.K.
}
\date{}

\maketitle
\noindent
\begin{abstract}

Non-random sample selection is a commonplace amongst many empirical studies and it appears when an output variable of interest is available only for a restricted non-random sub-sample of data. We introduce an extension of the generalized additive model which accounts for non-random sample selection by using a selection equation. The proposed approach allows for different distributions of the outcome variable, various dependence structures between the (outcome and selection) equations through the use of copulae, and nonparametric effects on the responses. Parameter estimation is carried out within a penalized likelihood and simultaneous equation framework. We establish asymptotic theory for the proposed penalized spline estimators, which extends the recent theoretical results for penalized splines in generalized additive models, such as those by \cite{Krivobokova-et-al-2009} and \cite{Yoshida-Naito-2014}. The empirical effectiveness of the approach is demonstrated through a simulation study.

\textbf{Key Words:} copula, generalized additive model, non-random sample selection, penalized regression spline, selection bias, simultaneous equation estimation.
\end{abstract}


\section{Introduction}

Non-random sample selection arises when an output variable of interest is available only for a restricted non-random sub-sample of the data. This often occurs in sociological, medical and economic studies where individuals systematically select themselves into (or out of) the sample (see, e.g., \cite{Guo-Fraser-2014}, chapter 4, \cite{Lennox-et-al-2012}, \cite{Vella-1998}, \cite{Collier-and-Mahoney-1996} and references therein). If the aim is to model an outcome of interest in the entire population and the link between its availability in the sub-sample and its observed values is through factors which can not be accounted for then any analysis based on the available sub-sample only will yield erroneous model structures and biased conclusions as the resulting inference may not extend to the unobserved group of observations. Sample selection models allow us to use the entire sample whether or not observations on the output variable were generated. In its classical form, it consists of two equations which model the probability of inclusion in the sample and the outcome variable through a set of available predictors, and of a joint bivariate distribution linking the two equations.

The sample selection model was first introduced by \cite{Gronau74}, \cite{Lewis74} and \cite{Heckman-1974}, in the context of estimating the number of working hours and wage rates of married women, some of whom did not work. \cite{Heckman76} formulated a unified approach to estimating this model as a simultaneous equation system. In the classical version, the error terms of the two equations are assumed to follow a bivariate normal distribution in which non-zero correlation indicates the presence of non-random sample selection. \cite{Heckman-1979} then translated the issue of sample selection into an omitted variable problem and proposed a simple and easy to implement estimation method known as two-step procedure. However, the method was proved to strongly rely on the assumption of normality and thus has been criticized for its lack of robustness to distributional misspecification and outliers \citep[e.g.,][]{Paarsch-1984, Little-and-Rubin-1987, Zuehlke-and-Zeman-1991,Zhelonkin-2013}. 

Various modifications and generalizations of the classical sample selection model have been proposed in the literature and we mention some of them. A non-parametric approach, which lifts the normality assumption, can be found in \cite{Das-Newey-Vella-2003}. Here a two-stage Heckman's method is adopted and a linear regression with the Heckman's selection correction is replaced with series expansions. Non-parametric methods are also considered in \cite{Lee08} and \cite{Chen10}. Semi-parametric techniques in a two-step scenario, similar to that in \cite{Das-Newey-Vella-2003}, are given in \cite{Ahn-and-Powell-1993}, \cite{Newey-1999} and \cite{Powell-1994}. Other semi-parametric approaches can be found in \cite{Gallant87}, \cite{Poweletal89}, \cite{Lee94a}, \cite{Lee94b}, \cite{Andrews98} and \cite{Newey09}. In the Bayesian framework, \cite{Chib09} deal with non-linear covariate effects using Markov chain Monte Carlo simulation techniques and simultaneous equation systems. \cite{bayesSampleSelection} further extend this approach by introducing a Bayesian algorithm based on low rank penalized B-splines for non-linear and varying-coefficient effects and Markov random-field priors for spatial effects. A frequentist counterpart of these Bayesian methods is discussed in \cite{Marra-and-Radice-2011} in the context of binary outcomes and \cite{Marra-and-Radice-2013} for the continuous outcome case. \cite{Zhelonkin-et-al-2012} introduce a procedure for robustifying the Heckman's two stage estimator by using M-estimators of Mallows' type for both stages. \cite{Marchenko-and-Genton-2012} and \cite{Ding2014} consider a bivariate Student-t distribution for the model's errors as a way of tackling heavy-tailed data. Several authors proposed specific copula functions to model the joint distribution of the selection and outcome equations; see, e.g., \cite{Prieger-2002} who employes an FGM bivariate copula or \cite{Lee-1983}. In the context of non-random sample selection, a more general copula approach, with a focus on Archimedean copulae, can be found in \cite{Smith-2003}. This stream of research is continued in \cite{Hasebe12} and \cite{Schwiebert13}. As emphasized by \cite{Genius}, copulae allow for the use of non-Gaussian distributions and has the additional benefit of making it possible to specify the marginal distributions independently of the dependence structure linking them. Importantly, while the copula approach is still fully parametric, it is typically computationally more feasible than non/semi-parametric approaches and it still allows the user to assess the sensitivity of results to different modeling assumptions.

The aim of this paper is to introduce a generalized additive model (GAM) which accounts for non-random sample selection. Thus, the classical GAM is extended by introducing an extra equation which models the selection process. The selection and outcome equations are linked by a joint probability distribution which is expressed in terms of a copula. Using different copulae allows us to capture different types of dependence while keeping the marginal distributions of the responses fixed. This approach is flexible and numerically tractable at the same time. Secondly, we make a step towards greater flexibility in modeling the relationship between predictors and responses by using a semiparametric approach in place of commonly used parametric formulae, thus capturing possibly complex non-linear relationships. We use penalized regression splines to approximate these relationships and employ a penalized likelihood and simultaneous equation estimation framework as presented in \cite{Marra-and-Radice-2013}, for instance. Penalized  splines were introduced in \cite{Sullivan-1986} and \cite{Eilers-Marx-1996} and practical aspects of their use received much attention ever since; see \cite{Rupp03} and \cite{Wood}. At the same time, the asymptotic theory for penalized spline estimators is relatively new and mainly focuses on models with one continuous predictor \citep[e.g.,][]{Hall-and-Opsomer-2005, Claeskens-et-al-2009, Wang-et-al-2011}. In the GAM context, the asymptotic normality of the penalized spline estimator has been proved by \cite{Krivobokova-et-al-2009} and recently extended by \cite{Yoshida-Naito-2014} to the case of any fixed number of predictors. This paper discusses the asymptotic properties of the penalized spline estimator of the proposed model. The asymptotic rate of the mean squared error of the linear predictor for the regression equation of interest is derived, which also allows us to derive its asymptotic bias and variance as well as its approximate distribution. The theoretical results are obtained for a general case of any fixed number of predictors and when the spline basis increases with the sample size. We show that even though the model structure is more complex than that of a classical GAM, the asymptotic properties of penalized spline estimators are still valid. Thus, the results extend the theoretical foundation of GAMs established so far. 

The paper is organized as follows. Section \ref{Sample-Selection-Models} briefly describes the classical sample selection model as well as Heckman's estimation procedure. Section \ref{Generalized-sample-selection-model} introduces a generalized sample selection model based on GAMs. Section \ref{Estimation} succinctly describes the estimation approach which is based on a penalized likelihood and simultaneous equation framework. In Section \ref{Asymptotics}, the asymptotic properties of the proposed penalized spline estimator are established. The finite-sample performance of the approach is investigated in Section \ref{Simulation-study}. Section \ref{discussion} discusses some extensions.


\section{Classical sample selection model}
\label{Sample-Selection-Models}

The classical sample selection model introduced by \cite{Heckman-1974} and \cite{Heckman-1979} is defined by the system of the following equations
\begin{eqnarray}
Y_{1i}^* = & {\bf x}_{i}^{(1)}\boldsymbol{\beta}_1 + \varepsilon_{1i} \label{1a}, \\
Y_{2i}^* = & {\bf x}_{i}^{(2)}\boldsymbol{\beta}_2 + \varepsilon_{2i} \label{1b}, \\
Y_{1i} = & \mathbbm{1}(Y_{1i}^*>0) \label{1d},\\
Y_{2i} = & Y_{2i}^* Y_{1i} \qquad \label{1c}, 
\end{eqnarray}
for $i=1,\ldots,n$, where $Y_{1i}^*$ and $Y_{2i}^*$  are unobserved latent variables and only values of $Y_{1i}$ and $Y_{2i}$ are observed. The symbol $\mathbbm{1}(\cdot)$ denotes throughout an indicator function.  Equation (\ref{1a}) is the so-called selection equation determining whether the observation of $Y_{2i}^*$ is generated, and equation (\ref{1b}) is the output equation of primary interest. Row vectors ${\bf x}^{(1)}=(x_{1}^{(1)},\ldots,x_{D_1}^{(1)})\in\R^{D_1}$ and ${\bf x}^{(2)}=(x_{1}^{(2)},\ldots,x_{D_2}^{(2)})\in\R^{D_2}$ contain predictors' values and are observed for the entire sample whereas $\boldsymbol{\beta}_1\in\R^p$ and $\boldsymbol{\beta}_2\in\R^q$ are unknown parameter vectors. The classical sample selection model assumes that the error terms $(\varepsilon_{1i},\varepsilon_{2i})$ follow the bivariate normal distribution
$$
\left[ \begin{array}{c} \varepsilon_{1i} \\ \varepsilon_{2i} \end{array} \right] \sim N \left(  \left[\begin{array}{cc} 0 \\ 0  \end{array} \right], 
\left[\begin{array}{cc} 1 & \rho\sigma \\ \rho\sigma & \sigma^2 \end{array}\right] \right),
$$
for $i=1,\ldots,n$. The variance of $\varepsilon_{1i}$ is assumed to be equal to 1 for the usual identification reason. The normality of $\varepsilon_{1i}$ implies that the selection equation represents a probit regression expressed in terms of a latent variable. The assumption of bivariate normality of error terms has been commonly used since. It implies that the log-likelihood of model (\ref{1a})-(\ref{1c}) can be expressed as (cf. \cite{Amemiya-1985} and eq. (10.7.3) therein)
\begin{eqnarray}
  \lefteqn{\ell(\boldsymbol{\beta}_1,\boldsymbol{\beta}_2,\sigma,\rho|Y_{11},\ldots,Y_{1n},Y_{21},\ldots,Y_{2n}) = 
	\sum_{i=1}^n \left\{ (1-Y_{1i}) \log \left(1-\Phi({\bf x}_{i}^{(1)}\boldsymbol{\beta}_1)\right) \right.} \nonumber\\
   &  & \left. + Y_{1i} \left[  
   \log\Phi\left( \frac{{\bf x}_{i}^{(1)}\boldsymbol{\beta}_1 + \frac{\rho}{\sigma}(Y_{2i}-{\bf x}_{i}^{(2)}\boldsymbol{\beta}_2)}{\sqrt{1-\rho^2}} \right) 
  + \log\phi\left( \frac{Y_{2i}-{\bf x}_{i}^{(2)}\boldsymbol{\beta}_2}{\sigma} \right) -\log\sigma  \right] \right\}, \nonumber
\end{eqnarray}
where $\Phi(\cdot)$ and $\phi(\cdot)$ denote the standard normal distribution function and density, respectively. Under the model assumptions, maximization of the above function, referred to as full-information maximum likelihood, results in consistent estimators of the parameters $\boldsymbol{\beta}_1$, $\boldsymbol{\beta}_2$, $\sigma$ and $\rho$ but until recently was computationally cumbersome. As a consequence \cite{Heckman-1979} proposed a two-step procedure, also known as limited-information maximum likelihood method. It is based on the fact that under the assumption of normality of the error terms the following equality holds
$$
  \E\left( \varepsilon_{2i} | Y_{1i}^*>0 \right) =  \sigma\rho \xi\left( {\bf x}_{i}^{(1)}\boldsymbol{\beta}_1 \right),
$$
where $\xi(u) = \frac{\phi(u)}{\Phi(u)}$ is the so-called inverse Mills ratio. Thus the conditional expectation of the outcome variable given that its value is observed equals
$$
  \E\left(Y_{2i}^* |  Y_{1i}^*>0  \right) = {\bf x}_{i}^{(2)}\boldsymbol{\beta}_2 + \sigma\rho \xi\left( {\bf x}_{i}^{(1)}\boldsymbol{\beta}_1 \right).
$$
In the first step of Heckman's procedure, probit model (\ref{1a}) is fitted in order to obtain an estimator of $\boldsymbol{\beta}_1$ and hence an estimator $\hat\xi_i=\xi\left( {\bf x}_{i}^{(1)}\hat{\boldsymbol{\beta}_1} \right)$ of $\xi\left( {\bf x}_{i}^{(1)}\boldsymbol{\beta}_1 \right)$. In the second step the following regression equation is estimated through the ordinary least squares method
$$
  \qquad Y_{2i} = {\bf x}_{i}^{(2)}\boldsymbol{\beta}_2 + \gamma \hat\xi_i + \tilde\varepsilon_{2i}, \quad i\in\{j: Y_{1j}^*>0 \},
$$
using the selected sub-sample only, where $\hat\xi_i$ is an added variable and $\gamma=\sigma\rho$ is the corresponding unknown parameter. After obtaining estimators $\hat{\boldsymbol{\beta}}_2$ and $\hat\gamma$, parameter $\sigma$ can be estimated as $\hat\sigma=\left( n_s^{-1} \sum {\hat{\tilde\varepsilon}}_{2i}^2 + n_s^{-1}\sum\hat\xi_i\left(\hat\xi_i+{\bf x}_{i}^{(1)}\boldsymbol{\beta}_1\right) \right)^{1/2}$ where $n_s$ is the size of the sub-sample, i.e. the cardinality of $\{j: Y_{1j}^*>0\}$ (see, for instance, \cite{Toomet-and-Henningsen-2008}). Thus the estimator of the correlation coefficient $\rho$ is defined as $\hat\rho=\hat\gamma/\hat\sigma$.

Although the inverse Mills ratio $\xi(u)$ is known to be nonlinear, in practice it often turns out to be approximately linear for most values in the $u$ range. This may lead to identification issues if both equations include the same set of predictors \citep[][p. 57]{Puhani-2000}. Thus, an exclusion restriction, which requires the availability of at least one predictor which is related to selection process but that has no direct effect on the outcome, has to be used in empirical applications. As mentioned in the introduction, Heckman's model received a number of criticisms due to its lack of robustness to departures from the model assumptions \citep[see, e.g.,][and references therein]{Paarsch-1984, Little-and-Rubin-1987,Zuehlke-and-Zeman-1991, Zhelonkin-2013}.


\section{Generalized additive sample selection model}
\label{Generalized-sample-selection-model}

We structure the proposed sample selection model in the following way. We first assume that the outcome variable of interest can be described by a generalized additive model \citep{Hast:Tibs:1990}. Then, to take the selection process into account we extend the model by adding a selection equation, which is also in the form of a generalized additive model. Finally, the two model equations are linked using a bivariate copula.


\subsection{Random component}
Let $F_{1}$ and $F_{2}$ denote the cumulative distribution functions of the latent selection variable $Y_{1}^*$ and output variable of interest $Y_{2}^*$, respectively. Analogically, $f_{1}$ and $f_{2}$ denote the density functions of $Y_{1}^*$ and $Y_{2}^*$.  We assume that $Y_{2}^*$ has density that belongs to an exponential family of distributions, i.e. it is of the form
\begin{equation}
	f_2(y_2|\eta_2,\phi)=\exp\left\{ \frac{y_2\eta_2-b(\eta_2)}{\phi} + c(y_2,\phi) \right\},
\label{expfam0}
\end{equation}
for some specific functions $b(\cdot)$ and $c(\cdot)$, where $\eta_2$ is the natural parameter and $\phi$ is the scale parameter. For simplicity, we assume that $\phi\equiv1$, so that
\begin{equation}
	f_2(y_2|\eta_2,\phi)=f_2(y_2|\eta_2)=\exp\left\{ y_2\eta_2-b(\eta_2) + c(y_2) \right\}.
\label{expfam}
\end{equation}
 It holds $\E(Y_2^*)=b'(\eta_2)$ and ${\rm Var}(Y_2^*)=b''(\eta_2)$, where $b'(\cdot)$ and $b''(\cdot)$ are the first and second derivatives of function $b(\cdot)$, respectively \citep[][p. 38]{vandervaart}.

Moreover, we assume that the latent selection variable $Y_1^*$ follows a normal distribution with mean $\eta_1$ and variance equal to $1$,
\begin{equation}
	f_1(y_1|\eta_1)= \exp\left( -(y_1-\eta_1)^2 \right).
\label{seld}
\end{equation}
The assumption ${\rm Var}(Y_1^*)=1$ is needed for the usual identification purpose. Then the observables are defined as before,
$$
\begin{array}{ll}
Y_{1} = & \mathbbm{1}(Y_{1}^*>0), \\
Y_{2} = & Y_{2}^* Y_{1}, \qquad 
\end{array}
$$
implying the probit regression model for the selection variable $Y_1$.

We specify the dependence structure between the two variables by taking advantage of the Sklar's theorem \citep{sklar}. It states that for any two random variables there exists a two-place function, called copula, which represents the joint cumulative distribution function of the pair in a manner which makes a clear distinction between the marginal distributions and the form of dependence between them. Thus a copula is a function that binds together the margins in order to form the joint cdf of the pair $(Y_1^*,Y_2^*)$. An exhaustive introduction to copula theory can be found in \cite{Nelsen-2006}, \cite{Joe-1997} and \cite{Schweizer-1991}. As in our model we would like to be able to specify the marginal distributions of $Y_1^*$ and $Y_2^*$, the use of copulae is a convenient approach which allows us to achieve modeling flexibility. Often copulae are parametrized with the so called association parameter $\theta$ which while varying leads to families of copulae with different strength of dependence. Thus, we use the symbol $C_\theta(\cdot,\cdot)$ throughout to denote a copula parametrized with $\theta$. Let the function $C_\theta$ be the copula such that the joint cdf of $(Y_{1}^*,Y_{2}^*)$ is equal to
\begin{equation}
  F(y_1,y_2) = C_\theta\left(F_{1}(y_1),F_{2}(y_2)\right).
\label{Cp}
\end{equation}
Function $C_\theta$ always exists and is unique for every $(y_1,y_2)$ in the support of the joint distribution $F$. Then the joint density of $(Y_{1}^*,Y_{2}^*)$ takes the form:
$$
  f(y_1,y_2) = \left. \frac{\partial ^2}{\partial u\partial v} C_{\theta}(u,v) \right|_{\tiny\substack{u=F_1(y_1)\\ v=F_2(y_2)}} f_1(y_1) f_2(y_2).
$$
The log-likelihood function for such defined sample selection model can be obtained by conditioning with respect to the value of the selection variable $Y_1$ (cf. \cite{Smith-2003}, p. 108). If $Y_1=0$ then the likelihood takes the simple form of $\mathbb{P}(Y_{1}=0)$, which is equivalent to $F_{1i}(0)$. Otherwise, the likelihood can be expressed, using the multiplication rule, as $\mathbb{P}(Y_{1}^*>0)f_{2|1}(y_{2}|Y_{1}^*>0)$, where $f_{2|1}$ denotes the conditional probability density function of $Y_{2}$ given $Y_{1}^*>0$. After substituting the conditional density $f_{2|1}(y_{2}|Y_{1}^*>0)$ by $\frac{1}{P(Y_{1}^*>0)} \frac{\partial}{\partial y_2}\left( F_{2}(y_2)-F(0,y_2) \right)$, we obtain the following log-likelihood
$$
  \ell =  (1-Y_{1}) \log F_{1}(0) + Y_{1} \log \left( f_{2}(Y_{2}) - \frac{\partial}{\partial y_2} F(0,y_2)\big|_{y_2\to Y_{2}} \right),
$$
from which, using (\ref{Cp}), we obtain 
$$
  \ell =  (1-Y_{1}) \log F_{1}(0) + Y_{1}  \log \left( f_{2}(Y_{2}) \left( 1 - z(Y_{2},\eta_{1},\eta_{2}) \right)  \right),
$$
where 
$$
  z(y_{2},\eta_{1},\eta_{2}) = \frac{\partial}{\partial v} C_{\theta}(F_{1}(0),v)\big|_{v\to F_{2}(y_{2})}.
$$
The function $z$ can be also expressed as 
$$
  z(y_{2},\eta_{1},\eta_{2}) = \mathbb{P}(Y_1^* \leq 0) \frac{f_{2|1}(y_2|Y_1^* \leq 0)}{ f_2(y_2)}.
$$
Thus it is directly related to the conditional distribution of the variable $Y_{2}^*$ for the unobserved data.

Using (\ref{expfam}), the log-likelihood can be written as
\begin{equation}
\label{loglik}
  \ell = (1-Y_{1}) \log F_{1}(0) + Y_{1} ( \eta_{2} Y_{2}  - b(\eta_{2}) +c(Y_{2}) + \log \left( 1 - z(Y_{2},\eta_{1},\eta_{2}) \right).
\end{equation}
The fact that $\E(Y_{2}) = b_2'(\eta_{2})$ implies
$$
\frac{\partial}{\partial \eta_2} \ell = Y_{1}(Y_{2}-\mu_{2}) + Y_{1} \frac{1}{1-z(Y_{2},\eta_{1},\eta_{2})} \frac{\partial}{\partial \eta_2} z(Y_{2},\eta_{1},\eta_{2}),
$$
where $\mu_{2}=\E(Y_{2})$.
Note that the first term in the expression above is equal to the score for a classical generalized linear model, i.e. when the sample selection does not appear and thus $Y_1$ always equals 1 and $z(Y_2,\eta_1,\eta_2)=0$. The second term corrects the score for sample selection bias. Interestingly, using the fact that the expected value of the score equals zero, we obtain that the expected value of the second term of the score equals minus the covariance between $Y_1$ and $Y_2$, i.e.
$$
 {\rm Cov}(Y_1,Y_2) = -\E\left (Y_1 \frac{\frac{\partial}{\partial \eta_2} z(Y_{2},\eta_{1},\eta_{2})}{1-z(Y_{2},\eta_{1},\eta_{2})} \right) 
$$
which provides another interesting interpretation of function $z(Y_2,\eta_1,\eta_2)$.

\subsubsection{Archimedean copulae}

\begin{figure}[ht!]
\begin{center}
\quad Normal \qquad\qquad\qquad\quad Clayton\\
\includegraphics[width=4.7cm]{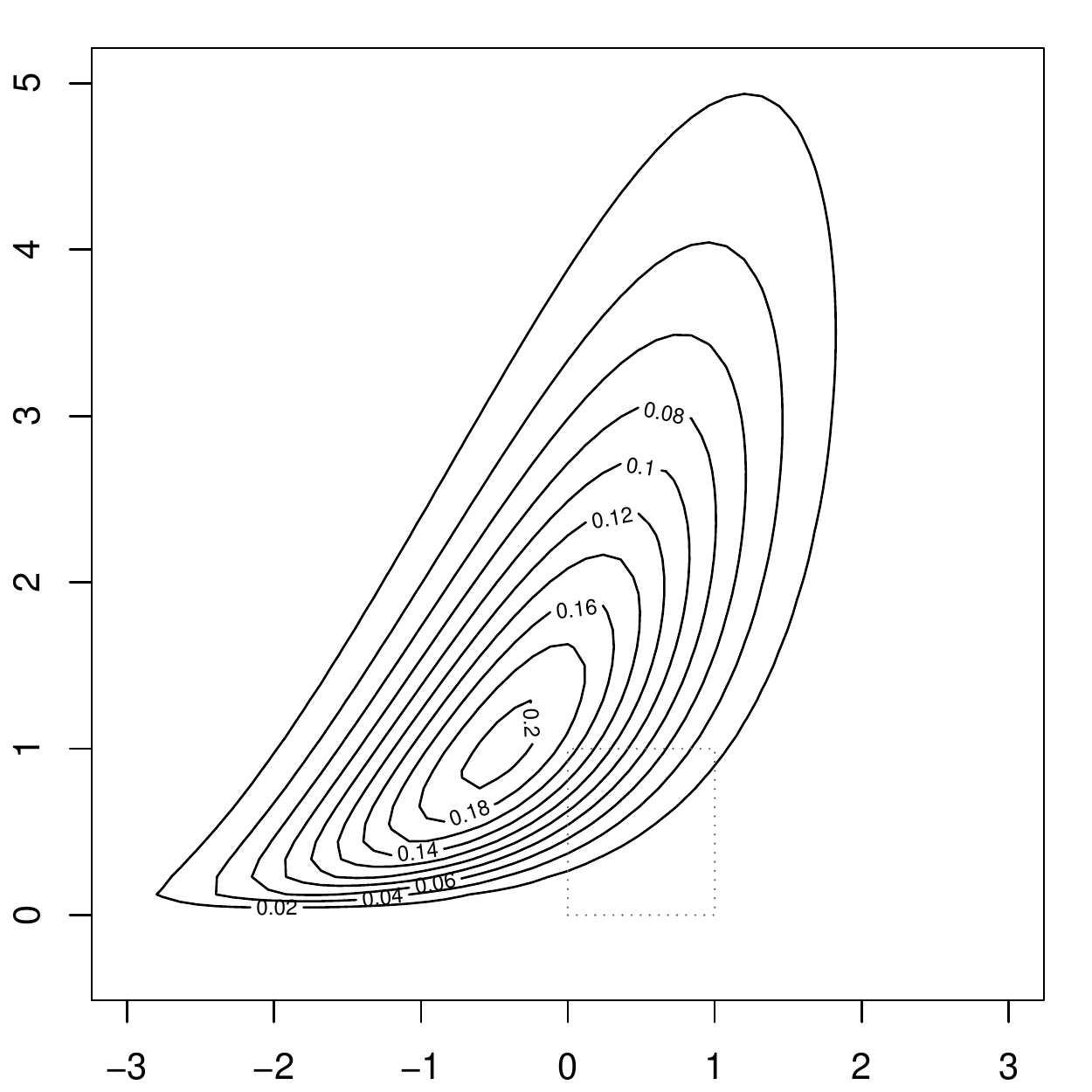}
\includegraphics[width=4.7cm]{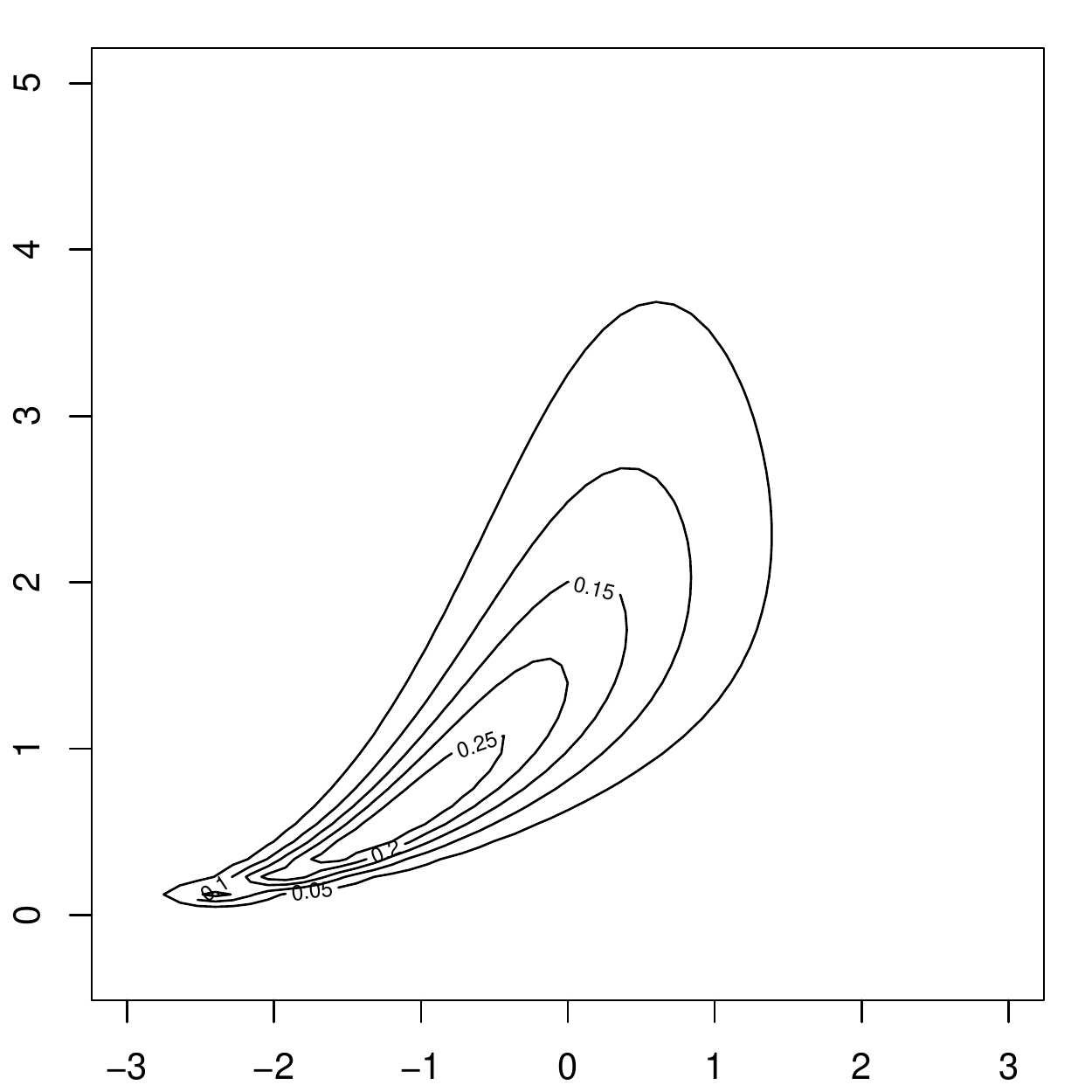}\\
\quad Joe \qquad\qquad\qquad\qquad Frank \\
\includegraphics[width=4.7cm]{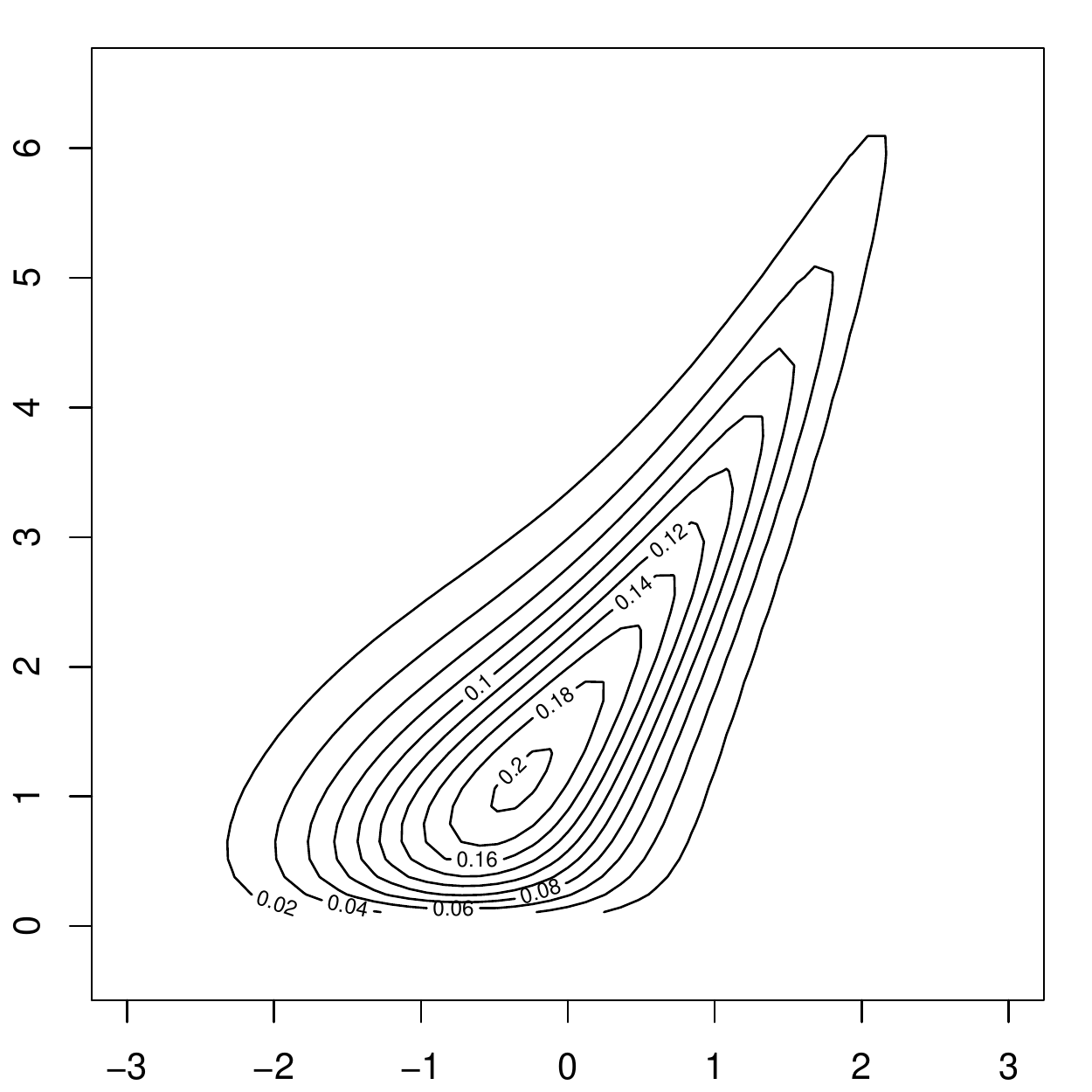}
\includegraphics[width=4.7cm]{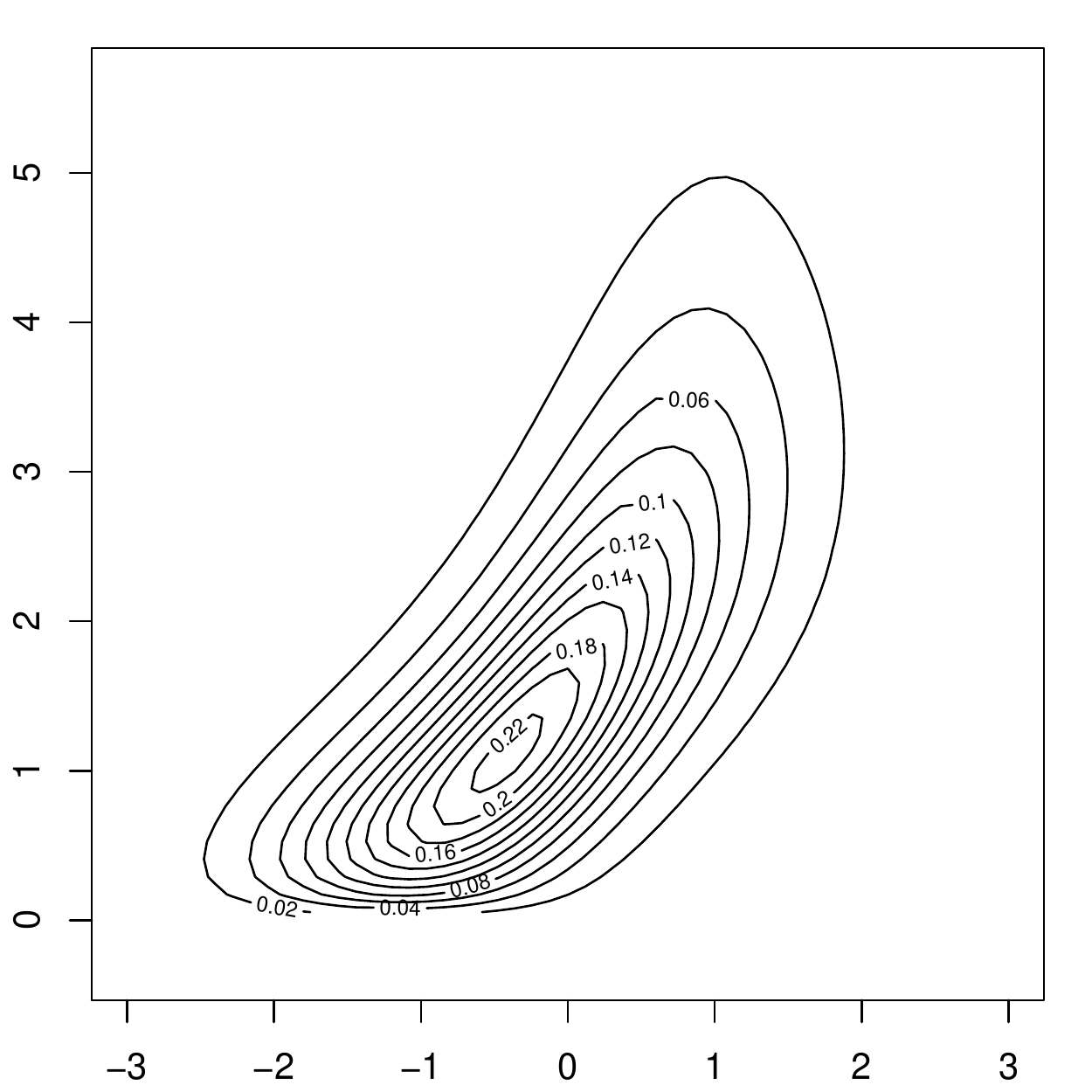}
\end{center}
\caption{Contour plots of bivariate densities for different copulae. The margins follow a standard normal distribution and a gamma distribution with shape 2 and rate 1. The value of dependence parameter $\theta$ is set to provide Kendall's $\tau=0.5$ for all four distributions.}
\label{contours}
\end{figure}

Although any copula can be used to link the two model equations (\ref{expfam}) and (\ref{seld}), the class of Archimedean copulae is particularly attractive for fitting the model in practice. This is because it provides many useful distributions that posses the advantage of analytical simplicity and dimension reduction because they are generated by a given function $\varphi:[0,1]\to[0,\infty)$ of only one argument such that
$$
  \varphi(C_{\theta}(u,v)) = \varphi(u)+\varphi(v).
$$
Function $\varphi$ is called a generator function and is assumed to be additive, continuous, convex, decreasing and meets the condition $\varphi(1)=0$. Table \ref{arcop} lists several popular bivariate Archimedean copulae whereas Figure \ref{contours} shows the contour plots of bivariate densities for normal and three Archimedean copulae (Clayton, Joe and Frank).
\\
For Archimedean copulae we have
$$
f_{2|1}(y_2|Y_1=1) = \frac{1}{\mathbb{P}(Y_1=1)} f_2(y_2) \left( 1-\frac{\varphi'(F_2)}{\varphi'(C_\theta)} \right),
$$
where $F_2=F_2(y_2)$ and $C_{\theta}=C_{\theta}(F_1(0),F_2(y_2))$.
Thus
$$
\E(Y_2^*|Y_1=1) = \frac{1}{\mathbb{P}(Y_1=1)} \left( \E\left(Y_2^*\right) - \int_{-\infty}^\infty yf_2(y)\frac{\varphi'(F_2)}{\varphi'(C_\theta)} dy\right).
$$
Hence the selection bias is equal to
$$
\frac{1}{\mathbb{P}(Y_1=1)} \left( \mathbb{P}(Y_1=0)\E \left(Y_2^*\right) - \int_{-\infty}^\infty yf_2(y)\frac{\varphi'(F_2)}{\varphi'(C_\theta)} dy\right).
$$
\vspace{2mm}

\begin{table}[ht!]
\begin{center}
\footnotesize
\begin{tabular}{lccc}
\hline
Name 	\qquad\qquad	& $C_{\theta}(u,v)$  & Parameter space & Generator $\varphi(t)$ \\\hline
Clayton		& $\left( u^{-\theta} + v^{-\theta} -1 \right)^{-1/\theta} $   & $\theta\in(0,\infty)$ & $\left(t^{-\theta}-1\right)\theta^{-1}$ \\
Joe				& $1-\left[ (1-u)^\theta + (1-v)^\theta - (1-u)^\theta(1-v)^\theta \right]^{1/\theta} $    & $\theta\in(1,\infty)$ & $-\log\left(1-(1-t)^\theta\right)$ \\
Frank     & $-\theta^{-1} \log \left[ 1+(e^{-\theta u}-1)(e^{-\theta v}-1)/(e^{-\theta}-1) \right] $  & $\theta\in\mathbb{R}\backslash \left\{0\right\}$ & $-\log\frac{e^{-\theta t}-1}{e^{-\theta}-1}$ \\
Gumbel		& $\exp\left\{-\left[(-\log u)^\theta + (-\log v)^\theta \right]^{1/\theta} \right\} $ & $\theta\in[1,\infty)$ & $(-\log t)^\theta$ \\
AMH				& $uv/\left[1-\theta(1-u)(1-v)\right]$ & $\theta\in[-1,1]$ & $\log\frac{1-\theta(1-t)}{t}$ \\\hline
\end{tabular}
\end{center}
\caption{Families of bivariate Archimedean copulae, with corresponding parameter range of the association parameter $\theta$ and generator $\varphi(t)$.}
\label{arcop}
\end{table}

\subsection{Systematic component}

Terms $\eta_1$ and $\eta_2$ are assumed to depend on sets of predictors ${\bf x}^{(1)}$ and ${\bf x}^{(2)}$, respectively, so that $\eta_1=\eta_1({\bf x}^{(1)})$ and $\eta_2=\eta_2({\bf x}^{(2)})$, where ${\bf x}^{(1)}=(x_{1}^{(1)},\ldots,x_{D_1}^{(1)})$ and ${\bf x}^{(2)}=(x_{1}^{(2)},\ldots,x_{D_2}^{(2)})$.
Moreover, we assume the following additive form of the functions $\eta_1({\bf x}^{(1)})$ and $\eta_2({\bf x}^{(2)})$
$$
\eta_1({\bf x}^{(1)}) = \eta_1^{(1)}(x_1^{(1)}) + \eta_2^{(1)}(x_2^{(1)}) + \ldots + \eta_{D_1}^{(1)}(x_{D_1}^{(1)}),
$$
and
$$
\eta_2({\bf x}^{(2)}) = \eta_1^{(2)}(x_1^{(2)}) + \eta_2^{(2)}(x_2^{(2)}) + \ldots + \eta_{D_2}^{(2)}(x_{D_2}^{(2)}).
$$
Functions $\eta_1({\bf x}^{(1)}) $ and $\eta_2({\bf x}^{(2)})$ are unknown.
We use spline basis functions to represent the unknown smooth functions. Specifically, we consider B-splines (\cite{De-Boor-2001}), which have attractive numerical and theoretical properties.
\\
\\
{\bf B-spline approximation}

\vspace{2mm}\noindent
On interval $[0,1]$ define sequence of knots $0=\kappa_0<\kappa_1<\ldots<\kappa_K=1$ and another $2p$ knots $\kappa_{K}=\kappa_{K+1}=\ldots=\kappa_{K+p}$ and $\kappa_{-p+1}=\kappa_{-p+2}=\ldots=\kappa_{-1}=\kappa_0$. Then B-spline basis functions of degree $p$ are defined recursively as
$$
B_{k,p}(x)= \frac{x-\kappa_{k-1}}{\kappa_{k+p-1}-\kappa_{k-1}} B_{k,p-1}(x) + 
            \frac{\kappa_{k+p}-x}{\kappa_{k+p}-\kappa_{k}} B_{k+1,p-1}(x)
$$
for $k=-p+1,\ldots,K$, with
$$
B_{k,0}(x) = \left\{ \begin{array}{ll} 1, & \kappa_{k-1}\leq x < \kappa_k\\
                                    0  & \text{otherwise}.  \end{array} \right.
$$
This gives $K+p$ basis functions $B_{-p+1,p}(x)$, \ldots, $B_{K,p}(x)$.
\\
\\
We approximate $\eta_j^{(1)}(x)$ by a linear combination of basis functions:
$$
	\sum_{k=-p+1}^K \alpha_{K,j} B_{k,p}(x) \quad\text{for } j=1,\ldots,D_1,
$$
and, similarly, $\eta_j^{(2)}(x)$ by
$$
	\sum_{k=-p+1}^K \beta_{K,j} B_{k,p}(x) \quad\text{for } j=1,\ldots,D_2.
$$
In the remaining of the paper we omit the subscript $p$ of basis functions $B_{k,p}$ and we use symbols $B_{-p+1}(x)$, \ldots, $B_{K}(x)$ to denote $p$-th B-spline basis functions. We define vectors $\boldsymbol{\alpha}_j \in\R^{p+K}$ and $\boldsymbol{\beta}_j  \in \R^{p+K}$ as
$$
\boldsymbol{\alpha}_j = (\alpha_{-p+1,j}, \ldots, \alpha_{K,j})^T \quad\text{for } j=1,\ldots,D_1
$$
and
$$
\boldsymbol{\beta}_j = (\beta_{-p+1,j}, \ldots, \beta_{K,j})^T \quad\text{for } j=1,\ldots,D_2
$$
and vectors 
$$
\boldsymbol{\alpha} = (\boldsymbol{\alpha}_1^T, \ldots, \boldsymbol{\alpha}_{D_1}^T)^T \in\R^{D_1(p+K)}
$$
and
$$
\boldsymbol{\beta} = (\boldsymbol{\beta}_1^T, \ldots, \boldsymbol{\beta}_{D_1}^T)^T \in\R^{D_2(p+K)}.
$$
Moreover, we use the following notation throughout. For a given $n\in\N$, assume that $(Y_{1i},Y_{2i})_{i=1}^n$ are independent random variables related to predictors' values ${\bf x}_i^{(1)} = (x_{1i}^{(1)}$, \ldots, $x_{D_1i}^{(1)})$ and ${\bf x}_i^{(2)} = ( x_{1i}^{(2)}$, \ldots, $x_{D_2i}^{(2)})$ for $i=1,\ldots,n$ such that $Y_{1i}=\mathbbm{1}(Y_{1i}^*>0)$ and $Y_{2i}=Y_{2i}^* Y_{1i}$, where $Y_{1i}^*$ has density (\ref{seld}) with $\eta_1=\eta_1({\bf x}_i^{(1)})$ and $Y_{2i}^*$ is distributed according to (\ref{expfam}) with $\eta_2=\eta_2({\bf x}_i^{(2)})$.

Let $F_{1i}$, $F_{2i}$ denote the distribution functions of $Y_{1i}^*$ and $Y_{2i}^*$ and let $F_i(\cdot,\cdot)$ be the joint cdf of the pair $(Y_{1i}^*,Y_{2i}^*)$. Moreover, for $j=1,\ldots,D_1$ let $X_j^{(1)}: n\times(p+K)$ be a matrix defined through
$$
X_j^{(1)} = \left[ B_{-p+k}(x_{ji}^{(1)}) \right]_{k=1,\ldots,K, \, i=1,\ldots,n}
$$
and for $j=1,\ldots,D_2$ let $X_j^{(2)}: n\times(p+K)$ be a matrix defined as
$$
X_j^{(2)} = \left[ B_{-p+k}(x_{ji}^{(2)}) \right]_{k=1,\ldots,K, \, i=1,\ldots,n}.
$$
Then let ${\bf X}^{(1)}: n\times D_1(p+K)$ and  ${\bf X}^{(2)}: n\times D_2(p+K)$ equal 
$$
 {\bf X}^{(1)} = \left[ X_1^{(1)},\ldots,  X_{D_1}^{(1)} \right]
$$
and
$$
 {\bf X}^{(2)} = \left[ X_1^{(2)},\ldots,  X_{D_2}^{(2)} \right].
$$
With the B-spline approximation, we postulate the parametric model
\begin{equation}
\begin{array}{l}
Y_{1i}^* \sim N({\bf X}_i^{(1)}{\boldsymbol\alpha},1)
\\
Y_{2i}^* \sim f_{2i}(y_2|{\bf X}_i^{(2)}) = \exp\left( y_2 {\bf X}_i^{(2)} {\boldsymbol\beta} - b({\bf X}_i^{(2)} {\boldsymbol\beta})  + c(y_2) \right)
\end{array}
\label{model2}
\end{equation}
where ${\bf X}_i^{(1)}$ and ${\bf X}_i^{(2)}$ denote the $i$-th rows of the matrices ${\bf X}^{(1)}$ and ${\bf X}^{(2)}$, respectively.


\section{Some estimation details} \label{Estimation}

In order to estimate parameters $\boldsymbol\alpha$, $\boldsymbol\beta$ and $\theta$, we employ a penalized likelihood approach which is common for regression spline models. Based on (\ref{loglik}) and (\ref{model2}), the log-likelihood given a random sample $(Y_{1i},Y_{2i})_{i=1}^n$ equals
\begin{equation}
\label{loglik2}
  \ell(\boldsymbol\alpha,\boldsymbol\beta,\theta) = \sum_{i=1}^n (1-Y_{1i}) \log F_{1i}(0) + Y_{1i} ( {\bf X}_i^{(2)} \boldsymbol\beta Y_{2i}  - b({\bf X}_i^{(2)} \boldsymbol\beta) +c(Y_{2i}) + \log \left( 1 - z(Y_{2i},\boldsymbol\alpha,\boldsymbol\beta,\theta) \right),
\end{equation}
where $ z(Y_{2i},\boldsymbol\alpha,\boldsymbol\beta,\theta)=\frac{\partial}{\partial v} C_{\theta}(F_{1i}(0),v)\big|_{v\to F_{2i}(Y_{2i})}$.
The normality of $Y_{1i}^*$ implies that $F_{1i}(0) = \Phi(-\eta_{1i})$, where $\Phi$ denotes the standard normal distribution function. The penalized log-likelihood equals
$$
	\ell_p(\boldsymbol\alpha,\boldsymbol\beta,\theta)= \ell(\boldsymbol\alpha,\boldsymbol\beta,\theta) -\frac{1}{2} \sum_{j=1}^{D_1} \lambda_j^{(1)} \alpha_j^T \Delta_m^T\Delta_m \alpha_j - \frac{1}{2}\sum_{j=1}^{D_2} \lambda_j^{(2)} \beta_j^T \Delta_m^T\Delta_m \beta_j =
$$
\begin{equation}
= \ell(\boldsymbol\alpha,\boldsymbol\beta,\theta) -\frac{1}{2} \boldsymbol\alpha^T Q_m^{(1)}(\lambda^{(1)}) \boldsymbol\alpha - \frac{1}{2} \boldsymbol\beta^T Q_m^{(2)}(\lambda^{(2)}) \boldsymbol\beta,
\label{Penloglik}
\end{equation}
where $\Delta_m: (K+p-m)\times(K+p)$ is the $m$-th difference matrix (Marx and Eilers, 1998), 
$Q_m^{(1)}(\lambda^{(1)})={\rm diag}(\lambda_1^{(1)}\Delta_m^T\Delta_m, \ldots, \lambda_{D_1}^{(1)}\Delta_m^T\Delta_m)$,
$Q_m^{(2)}(\lambda^{(2)})={\rm diag}(\lambda_1^{(2)}\Delta_m^T\Delta_m, \ldots, \lambda_{D_2}^{(2)}\Delta_m^T\Delta_m)$, and $\lambda^{(1)}=(\lambda_1^{(1)},\ldots,\lambda_{D_1}^{(1)})$ and $\lambda^{(2)}=(\lambda_1^{(2)},\ldots,\lambda_{D_2}^{(2)})$ are smoothing parameters controlling the trade-off between smoothness and fitness. 

For a given $\bm\lambda = (\lambda^{(1)},\lambda^{(2)})$, we seek to maximize (\ref{Penloglik}). In practice, an iterative procedure based on a trust region approach can be used to achieve this. This method is generally more stable and faster than its line-search counterparts (such as Newton-Raphson), particularly for functions that are, for example, non-concave and/or exhibit regions that are close to flat; see \citet[][Chapter 4]{Nocedal} for full details. Such functions can occur relatively frequently in bivariate models, often leading to convergence failures \citep{Andrews-1999,Butler-1996,Chiburis-et-al-2012}. Some details related to the process of iterative maximization of the penalized likelihood via a trust region algorithm are presented in subsection \ref{trust} of Appendix \ref{Estimation_details}. 

Data-driven and automatic smoothing parameter estimation is pivotal for practical modeling, especially when the data are partly censored as in our case, and each model equation contains more than one smooth component. An automatic approach allows us to determine the shape of the smooth functions from the data, hence avoiding arbitrary decisions by the researcher as to the relevant functional form for continuous variables. For single equation spline models, there are a number of methods for automatically estimating smoothing parameters within a penalized likelihood framework; see \cite{Rupp03} and \cite{Wood} for excellent detailed overviews. In our context, we propose to use the smoothing approach based on Un-Biased Risk Estimator as detailed in subsection \ref{estimating_lambda} of Appendix \ref{Estimation_details}.


\section{Asymptotic theory} \label{Asymptotics}

In this section, the asymptotic consistency of the linear predictor $\hat\eta_2$ based on the penalized maximum likelihood estimators for the regression equation of interest is shown. In particular, the asymptotic rate of the mean squared error of $\hat\eta_2$ is derived. The theoretical considerations also allows to derive its asymptotic bias and variance and its approximate distribution.
\\
\\
First, we introduce the notation that we use throughout. Denote $\bdelta=(\balpha,\bbeta,\theta)$ and let
$$
  G_{n}(\bdelta) = \left( G_{n}^{\balpha}(\bdelta),G_{n}^{\bbeta}(\bdelta),G_{n}^{\theta}(\bdelta)   \right),
$$
where
$$
  G_{n}^{\balpha}(\bdelta)=\frac{\partial\ell}{\partial\balpha} (\bdelta)= \left( \frac{\partial\ell}{\partial\balpha_1},\ldots, \frac{\partial\ell}{\partial\balpha_{D_1}} \right) \in \R^{D_1(K+p+1)},
$$
$$
  G_{n}^{\bbeta}(\bdelta)=\frac{\partial\ell}{\partial\bbeta} (\bdelta)= \left( \frac{\partial\ell}{\partial\bbeta_1},\ldots, \frac{\partial\ell}{\partial\bbeta_{D_2}} \right)  \in \R^{D_1(K+p+1)},
$$
and
$$
  G_{n}^{\theta}(\bdelta) = \frac{\partial\ell}{\partial\theta}(\bdelta) \in\R.
$$
Analogically, let
$$
  G_{n,p}(\bdelta) = \left( G_{n,p}^{\balpha}(\bdelta),G_{n,p}^{\bbeta}(\bdelta),G_{n,p}^{\theta}(\bdelta) \right) = 
	\left(  \frac{\partial\ell_p}{\partial\balpha} (\bdelta),  \frac{\partial\ell_p}{\partial\bbeta} (\bdelta),   \frac{\partial\ell_p}{\partial\theta}(\bdelta)\right).
$$

\vspace{2mm}\noindent
Proceeding to the hessian matrix, let us denote
$$
  H_n^{\balpha}(\bdelta) = \frac{\partial^2 \ell}{\partial\balpha\partial\balpha^T} (\bdelta), \,\,
  H_n^{\bbeta}(\bdelta) = \frac{\partial^2 \ell}{\partial\bbeta\partial\bbeta^T} (\bdelta), \,\,
  H_n^{\bbeta,\balpha}(\bdelta) = \frac{\partial^2 \ell}{\partial\balpha\partial\bbeta^T} (\bdelta) \,\,
	\text{and so on.}
$$
Moreover, let
$$
  H_n(\bdelta) = \left[ \begin{array}{ccc}
	  H_n^{\balpha}(\bdelta)   &  H_n^{\balpha,\bbeta}(\bdelta) & H_n^{\balpha,\theta}(\bdelta) 
		\\
		H_n^{\bbeta,\balpha}(\bdelta)  & H_n^{\bbeta}(\bdelta) & H_n^{\bbeta,\theta}(\bdelta) 
		\\
		H_n^{\theta,\balpha}(\bdelta)  &  H_n^{\theta,\bbeta}(\bdelta)  & H_n^{\theta,\theta}(\bdelta) 
	\end{array} \right].
$$
Analogically, let 
$$
  H_{n,p}^{\balpha}(\bdelta) = \frac{\partial^2 \ell_p}{\partial\balpha\partial\balpha^T} (\bdelta) 
	= H_n^{\balpha}(\bdelta) - Q_m^{(1)}(\lambda_n^{(1)}),
$$
$$	
  H_{n,p}^{\bbeta}(\bdelta) = \frac{\partial^2 \ell_p}{\partial\bbeta\partial\bbeta^T} (\bdelta)
	= H_n^{\bbeta}(\bdelta) - Q_m^{(2)}(\lambda_n^{(2)}), 
$$
and so on, and $H_{n,p}(\bdelta) = \frac{\partial^2 \ell_p}{\partial\bdelta\partial\bdelta^T} (\bdelta)$.

\vspace{2mm}\noindent
Let
$$
  F_{n}^{\balpha}(\bdelta)=\E\left[ H_n^{\balpha}(\bdelta)  \right], \,\,
  F_{n}^{\bbeta}(\bdelta)=\E\left[ H_n^{\bbeta}(\bdelta)  \right], \,\,
	F_{n}^{\balpha,\bbeta}(\bdelta)=\E\left[ H_n^{\balpha,\bbeta}(\bdelta),  \right]
$$
$$
 F_{n,p}^{\balpha}(\bdelta)=F_{n}^{\balpha}(\bdelta) - Q_m^{(1)}(\lambda_n^{(1)}), \,\,
 F_{n,p}^{\bbeta}(\bdelta)=F_{n}^{\bbeta}(\bdelta) - Q_m^{(2)}(\lambda_n^{(2)}).
$$

\noindent Let $\bdelta^0$ denote a parameter vector that satisfies the condition $\E G_n(\bdelta^0) = 0$. Then $\bdelta^0$ is maximizer of the expected unpenalized log-likelihood and provides the best approximation of $(\eta_1,\eta_2)$ in terms of Kullback-Leibler measure as it minimizes Kullback-Leibler distance. 
We adopt the following assumptions:
\begin{description}
\item[A1.] All partial derivatives up to the order 3 of copula function $C_{\theta}(u,v)$ w.r.t  $u$, $v$ and $\theta$ exist and are bounded.
\item[A2.] The function $z(y_2,\balpha,\bbeta,\theta)$ is bounded away from $1$.
\item[A3.] $\max_{l=1,2;\, j=1,\ldots,D_l}\left(\lambda_j^{(l)} \right)= O\left( n^{\gamma}\right)$ where $\gamma\leq\frac{2}{2p+3}$.
\item[A4.] The explanatory variables ${\bf x}^{(1)}$ and ${\bf x}^{(1)}$ are distributed on unit cubes $[0,1]^{D_1}$ and $[0,1]^{D_2}$, respectively.
\item[A5.] The knots of the B-spline basis are equidistantly located so that $\kappa_k-\kappa_{k-1} = K_n^{-1}$ for $k=1,\ldots,K_n$ and the dimension of the spline basis satisfies $K_n=O(n^{1/(2p+3)})$.
\item[A6.] $K_n$ is such that $(D_1+D_2)(K_n+p)<n$.
\end{description}

\noindent
{\bf Theorem.} {\it Under assumptions (A1)-(A6) the estimate $\hat\eta(x)$ has asymptotic expansion
$$
 \hat\eta(x)-\eta^0(x) \approx {\bf X}(x) F_p^{-1}(\bdelta^0)G_p(\bdelta^0),
$$
which implies
$$
  {\rm MSE}(\hat\eta(x)) = \E(\hat\eta(x)-\eta^0(x))^2 = O(n^{-(2p+2)/(2p+3)}).
$$
}

\noindent 
Before we proceed to the proof of the theorem we derive analytic formulae for the gradient and hessian of the penalized log-likelihood as their properties will play a central role in the asymptotic considerations.
\\
\\
{\bf Gradient of the penalized likelihood}

\vspace{2mm}\noindent
Straightforward calculations yield
$$
  G_{n,p}^{\balpha}(\bdelta)= -\sum_{i=1}^n \left\{ (1-Y_{1i})\left( \Phi(-{\bf X}_i^{(1)}\balpha) \right)^{-1} - Y_{1i} \frac{\nabla C_i}{1-z_i}  \right\}
	\varphi(-{\bf X}_i^{(1)}\balpha) {\bf X}_i^{(1)} - \balpha^T Q_m^{(1)}(\lambda_n^{(1)}),
$$
where $\varphi(\cdot)$ is the density function of the standard normal distribution and \\ $\nabla C_i = \frac{\partial^2}{\partial u\partial v}C_\theta(u,v)|_{u=F_{1i}(0), v=F_{2i}(Y_{2i})}$ and $z_i=z(Y_{2i},\balpha,\bbeta,\theta)$. 
Moreover,
$$
  G_{n,p}^{\bbeta}(\bdelta) = \sum_{i=1}^n Y_{1i}\left[ Y_{2i}-b'({\bf X}_i^{(2)}\bbeta) - \frac{z_i'}{1-z_i} \right] {\bf X}_i^{(2)} - \bbeta^T Q_m^{(2)}(\lambda_n^{(2)}),
$$
where 
$$
  z_i' = \frac{\partial z_i}{\partial\eta_2} = \frac{\partial^2}{\partial v^2} C_{\theta}(F_{1i}(0),v)|_{v=F_{2i}(Y_{2i})} \frac{\partial F_{2i}}{\partial\eta_2},
$$
and
$$
  G_{n,p}^{\theta}(\bdelta) = -\sum_{i=1}^n Y_{1i} \frac{1}{1-z_i}\frac{\partial z_i}{\partial\theta},
$$
where $\frac{\partial z_i}{\partial\theta} = \frac{\partial^2}{\partial\theta\partial v} C_{\theta}(F_{1i}(0),v)|_{v=F_{2i}(Y_{2i})}$.
Thus we can write the gradients in a matrix form
$$
  G_{n,p}^{\balpha}(\bdelta)= {\bf a}^T {\bf X}^{(1)} - \balpha^T Q_m^{(1)}(\lambda_n^{(1)}),
$$
where ${\bf a} = (a_1,\ldots,a_n)$ with $a_i = -\left\{ (1-Y_{1i})\left( \Phi(-{\bf X}_i^{(1)}\balpha) \right)^{-1} - Y_{1i} \frac{\nabla C_i}{1-z_i}  \right\} \varphi(-{\bf X}_i^{(1)}\balpha)$ for $i=1,\ldots,n$. Analogically,
$$
  G_{n,p}^{\bbeta}(\bdelta) = {\bf b}^T {\bf X}^{(2)} - \bbeta^T Q_m^{(2)}(\lambda_n^{(2)}),
$$
where ${\bf b}=(b_1,\ldots,b_n)$ with $b_i=Y_{1i}\left[ Y_{2i}-b'({\bf X}_i^{(2)}\bbeta) - \frac{z_i'}{1-z_i} \right]$ for $i=1,\ldots,n$. Moreover,
$$
  G_{n,p}^{\theta}(\bdelta) = {\bf c}^T {\bf 1},
$$
where ${\bf c}=(c_1,\ldots,c_n)$ with $c_i=-Y_{1i} \frac{1}{1-z_i}\frac{\partial z_i}{\partial\theta}$ and ${\bf 1}=(1,\ldots,1)^T\in\R^n$.
\\
\\

\noindent
{\bf Hessian}

\vspace{2mm}\noindent
Straightforward but tedious calculations yield the following lemma.

\vspace{2mm}\noindent
{\bf Lemma 1.} The component matrices of $H_n(\bdelta)$ take the following forms:
$$
  H_n^{\balpha}(\bdelta) = \left({\bf X}^{(1)}\right)^T W_1 {\bf X}^{(1)}, \quad
	H_n^{\bbeta}(\bdelta) = \left({\bf X}^{(2)}\right)^T W_2 {\bf X}^{(2)}, \quad
	H_n^{\balpha,\bbeta}(\bdelta) = \left({\bf X}^{(1)}\right)^T W_3 {\bf X}^{(2)},
$$
$$
  H_n^{\theta,\balpha}(\bdelta)={\bf 1}^T W_4 {\bf X}^{(1)},  \quad
	H_n^{\theta,\bbeta}(\bdelta)={\bf 1}^T W_5 {\bf X}^{(2)},   \quad
	H_n^{\theta,\theta}(\bdelta) ={\bf 1}^T W_6 {\bf 1},
$$
where $W_j={\rm diag}(w_1^{(j)},\ldots,w_n^{(j)})$ for $j=1,\ldots,6$ and
\begin{equation}
  w_i^{(1)} = (1-Y_{1i})(F_{1i}(0))^{-1} (F_{1i}(0)^{-1} \varphi(-{\bf X}_i^{(1)}\balpha) -1) \varphi(-{\bf X}_i^{(1)}\balpha)+
\label{w1i}
\end{equation}
$$
  \qquad\qquad + Y_{1i}\left\{ \left[ \frac{1}{1-z_i} \frac{\partial^3 C_i}{\partial u^2\partial v} - \left( \frac{\nabla C_i}{1-z_i} \right)^2\right] \varphi(-{\bf X}_i^{(1)}\balpha) - \frac{\nabla C_i}{1-z_i}
	   \right\} \varphi(-{\bf X}_i^{(1)}\balpha),
$$
\begin{equation}
  w_i^{(2)} = Y_{1i} \left[ -b''({\bf X}_i^{(2)}\bbeta)  +  \frac{1}{(1-z_i)^2} \left( z_i''(1-z_i) + (z_i')^2 \right)\right],
\label{w2i}
\end{equation}
\begin{equation}
  w_i^{(3)} = Y_{1i}\left[ \frac{1}{1-z_i} \frac{\partial^3 C_i}{\partial u^2\partial v} + \frac{z_i'}{(1-z_i)^2} \nabla C_i \right] \,\,\text{where}\,\, z_i''=\frac{\partial z_i'}{\partial \eta_2},
\label{w3i}
\end{equation}
\begin{equation}
  w_i^{(4)} = \frac{\partial}{\partial\theta}\frac{\nabla C_i}{1-z_i}\phi(-{\bf X}_i^{(1)}\balpha),
\label{w4i}
\end{equation}
\begin{equation}
  w_i^{(5)} = - Y_{1i}\frac{\partial}{\partial\theta}\frac{z_i'}{1-z_i},
\label{w5i}
\end{equation}
\begin{equation}
  w_i^{(6)} = - Y_{1i}\frac{\partial}{\partial\theta}\left( \frac{1}{1-z_i} \frac{\partial z_i}{\partial\theta}\right).
\label{w6i}
\end{equation}
	
\noindent
{\bf Corollary 1.} The Hessian $H_{n}(\bdelta)$ can be written in the following matrix form
$$
  H_{n}(\bdelta)= 
	  \left[ \begin{array}{ccc} \left({\bf X}^{(1)}\right)^T & {\bf 0} & {\bf 0} \\ {\bf 0} &  \left({\bf X}^{(2)}\right)^T & {\bf 0} \\
		                          {\bf 0} & {\bf 0} & {\bf 1}^T   \end{array} \right]
	  \left[ \begin{array}{ccc} W_1 & W_3 &  W_4 \\ W_3 &  W_2 & W_5 \\ W_4 & W_5 & W_6 \end{array} \right]
	  \left[ \begin{array}{ccc} {\bf X}^{(1)} & {\bf 0} & {\bf 0} \\ {\bf 0} &  {\bf X}^{(2)} & {\bf 0} \\
		                          {\bf 0} & {\bf 0} & {\bf 1} \end{array} \right]
$$
where the matrices $W_j={\rm diag}(w_1^{(j)},\ldots,w_n^{(j)})$ for $j=1,\ldots,6$ are given by the expressions (\ref{w1i})-(\ref{w6i}).
\\
\\
In order to proof Theorem, we use several lemmas, stated below. 
\\
\\
{\bf Lemma 2.} Under assumptions (A1)-(A6), elements of the matrix $F_{n}(\bdelta^0)$ are of order $O\left(\frac{n}{K_n}\right)$ and elements of the matrix $F_{n,p}(\bdelta^0)$ are of order $O\left(\frac{n}{K_n} + \max_{l=1,2;\,j=1,\ldots,D_l}\lambda_j^{(l)} K_n^{2p}\right)$.
\\
\\

\noindent
{\bf Lemma 3.} Under assumptions (A1)-(A6), elements of the matrix $\left(F_{n,p}(\bdelta^0)\right)^{-1}$ are of order $O_P\left(\frac{K_n}{n}\right)$. 
\\
\\

\noindent
{\bf Lemma 4.} Elements of the matrix $H_{n}(\bdelta^0) - F_{n}(\bdelta_0)$ are of order $O_P\left( \frac{n}{k}\right)$.
\\
\\

\noindent
The proofs of Lemmas 2 - 4 are presented in Appendix \ref{proofs}.
\\
\\

\noindent
{\bf Proof of Theorem.}
First, we expand $G_p(\cdot)$ around $\bdelta^0$. Let $M_n=(D_1+D_2)(K_n+p)$. For $j=1,\ldots,D_1+D_2+1$
$$
  0 = \frac{\partial l_p}{\partial \delta^j}(\hat\bdelta) = \frac{\partial l_p}{\partial \delta^j}(\bdelta^0)+
	\sum_{l=1}^{M_n} \frac{\partial^2 l_p}{\partial \delta^j\partial\delta^l}(\bdelta^0)(\hat\bdelta^l - \bdelta^{0,l}) +
$$
$$
  + \sum_{l=1}^{M_n} \sum_{r=1}^{M_n} (\hat\bdelta^l - \bdelta^{0,l}) \frac{\partial^3 l_p}{\partial \delta^j\partial\delta^l\partial\delta^r}(\bdelta^0)(\hat\bdelta^r - \bdelta^{0,r}) + o( R_n),
$$
where $R_n = \sum_{l=1}^{M_n} \sum_{r=1}^{(D_1+D_2)(K_n+p)} (\hat\bdelta^l - \bdelta^{0,l}) \frac{\partial^3 l_p}{\partial \delta^j\partial\delta^l\partial\delta^r}(\bdelta^0)(\hat\bdelta^r - \bdelta^{0,r})$.
\\
Series inversion yields
$$
  \hat\delta^j - \delta_0^j = -\sum_{l=1}^{M_n} a_{jl} \frac{\partial l_p}{\partial \delta^l}(\bdelta^0) 
	+ \frac{1}{2} \sum_{l=1}^{M_n} \sum_{r=1}^{M_n}
	   b_{jlr} \frac{\partial l_p}{\partial \delta^l}(\bdelta^0) \frac{\partial l_p}{\partial \delta^r}(\bdelta^0) +\ldots
$$
where $a_{jl}$ is $(j,l)$-element of the inverse of matrix $H_p(\bdelta^0)$ and \\
$b_{jlr} = \sum_{s=1}^{M_n} \sum_{t=1}^{M_n} \sum_{u=1}^{M_n} a_{js}a_{lt}a_{ru} \frac{\partial^3 l_p}{\partial \delta^s\partial\delta^t\partial\delta^u}(\bdelta^0)$.
\\
Then 
$$
  \left( H_p(\bdelta^0) \right)^{-1} = \left( F_p(\bdelta^0) + \left( H_p(\bdelta^0) - F_p(\bdelta^0) \right) \right)^{-1} = 
	\left( F_p(\bdelta^0) + S(\bdelta^0) \right)^{-1} =
$$
$$
  = F_p(\bdelta^0)^{-1} - F_p(\bdelta^0)^{-1} S(\bdelta^0) F_p(\bdelta^0)^{-1} + F_p(\bdelta^0)^{-1} S(\bdelta^0) F_p(\bdelta^0)^{-1}S(\bdelta^0) F_p(\bdelta^0)^{-1} + \ldots=
$$
$$
  = F_p(\bdelta^0)^{-1} \left( I - S(\bdelta^0) F_p(\bdelta^0)^{-1} + S(\bdelta^0) F_p(\bdelta^0)^{-1}S(\bdelta^0) F_p(\bdelta^0)^{-1} + \ldots \right) = 
$$
$$
  = F_p(\bdelta^0)^{-1} \left( I+ O\left( \sqrt{K_n/n} \right) \right).
$$
Moreover, 
$$
  \frac{\partial^3 l_p}{\partial \delta^s\partial\delta^t\partial\delta^u}(\bdelta^0) = 
	\sum_{i=1}^n \frac{\partial w_i}{\partial \delta^s} B_{-p+t}(x_i) B_{-p+u}(x_i) = O_P(n/K_n),
$$
as $\sum_{i=1}^n B_{-p+t}(x_i) B_{-p+u}(x_i) = O(n/K_n)$ and $\frac{\partial w_i}{\partial \delta^s}$ is bounded. 
This yields $b_{jlr} = O(K_n^2/n^2)$ and in consequence
$$
 \sum_{l=1}^{M_n} \sum_{r=1}^{M_n}
	   b_{jlr} \frac{\partial l_p}{\partial \delta^l}(\bdelta^0) \frac{\partial l_p}{\partial \delta^r}(\bdelta^0) = o_p(K_n/n).
$$
Thus
$$
  \hat\bdelta^j - \bdelta^{0,j} = \sum_{l=1}^{M_n} a_{jl} \left(-\frac{\partial l_p}{\partial \delta^l}(\bdelta^0) \right) + o_P(K_n/n) = 
$$
$$
 = \sum_{l=1}^{M_n} \bar f_{jl} \left(-\frac{\partial l_p}{\partial \delta^l}(\bdelta^0) \right)(1+o(1)) + o_P(K_n/n),
$$
where $\bar f_{jl}$ is $(j,l)$-element of the matrix $F_p(\bdelta^0)^{-1}$. Hence we can write the above equation in a matrix form
\begin{equation}
  \hat\bdelta - \bdelta^{0} = - F_p(\bdelta^0)^{-1} G_p(\bdelta^0)(1+o(1)) + o_P(K_n/n).
\label{last}
\end{equation}
Thus assumptions (A3) and (A5) together with the fact that $\E(G(\bdelta^0))=0$ yield
$$
  \E\left|\left|\hat\eta(x)-\hat\eta^0(x)\right|\right|^2 = \E({\bf X}(x)\hat\bdelta - {\bf X}(x)\bdelta^0)({\bf X}(x)\hat\bdelta - {\bf X}(x)\bdelta^0)^T= O(K_n^{-(p+1)}+\lambda_{jn}K_n^{1-m}/n).
$$
Hence $\E||\hat\eta(x)-\hat\eta^0(x)||^2  = O(n^{-(2p+2)/(2p+3)})$.
\begin{flushright}
\vspace{-5mm}
$\square$
\end{flushright}

\noindent
{\bf Remark} Expansion (\ref{last}) yields the approximate variance of $\hat\eta(x)$ in the form
$$
  {\rm Var}(\hat\eta(x)) \approx {\bf X}(x)F_p(\bdelta^0)^{-1}F(\bdelta^0) F_p(\bdelta^0)^{-1}{\bf X}(x)^T,
$$
which can be shown to be of order $O(K_n/n)$ as $n\to \infty$. Moreover, using the Central Limit Theorem we obtain the approximate  distribution of $\hat\eta(x)$ as
$$
  \hat\eta(x) \stackrel{a}{\sim} N(\eta^0(x),{\rm Var}(\hat\eta(x)) ).
$$


\section{Simulations}
\label{Simulation-study}

In this section, the properties of the proposed generalized additive sample selection model are investigated empirically. Specifically, we first assess the effectiveness of the proposed approach at finite sample sizes and then provide some evidence of the potential inaccuracy arising from modeling transformed outcomes.  
 
\subsection{Empirical consistency}

Data were generated as follows. For the latent selection variable, it was assumed that
$$
Y_{1i}^*=\alpha_{0} + s_1(x_1) + s_2(x_2) + \alpha_{4} x_4 + \alpha_5 x_5 + \varepsilon_i,
$$
with $\alpha_{0}=0.7$, $s_{1}(x)=-0.2\sin(\frac{\pi}{46}x)$, $s_2(x)=-0.0004(x+0.01x^{1/3})$, $\alpha_4=0.6$, $\alpha_5 = -0.4$ and $\varepsilon_i\sim N(0,1)$, whereas the outcome variable $Y_{2i}^*$ was assumed to pertain to a gamma distribution with shape parameter $k=2$ and expected value $\mu_i=\E( Y_{2i}^*)$ such that 
$$
 \log( \mu_i ) = \beta_0 + s_3(x_1) + s_4(x_3) + \beta_4 x_4  + \beta_5 x_5,
$$
with $\beta_0=-1.5$, $s_3(x) = 0.0006 \exp(0.1x)$, $s_4(x) = 0.03x$, $\beta_4=-1$ and $\beta_5 = 0.75$. The two equations were linked using a Gumbel copula with association parameter $\theta=3$. The covariates were generated as $x_1\sim {\rm Uniform}(16,66)$, $x_2\sim {\rm Uniform}(10,70)$, $x_3\sim {\rm Uniform}(0,20)$ and $x_4$ and $x_5$ were binary variables taking values $0$ and $1$ with equal probabilities. 

To investigate the asymptotic behavior of the estimators, data sets of increasing sizes were considered: $n=500, \, 1000, \, 1500, \, 2000, \, 2500, \, 3000$. For each generated data set a generalized additive sample selection model assuming a gamma distribution for the outcome with log link function was fitted using the {\tt SemiParSampleSel} function from the package {\tt SemiParSampleSel} \citep{Wojtys-et-al-2015} in {\tt R} \citep{R-2014}. Additionally, a univariate generalized additive model based only on the observed outcomes was fitted for comparison using the {\tt gam} function from {\tt mgcv} \citep{Wood}. The two above models will be referred to as GASSM and GAM, respectively. The number of Monte Carlo repetitions was $200$. 

\begin{figure}
\begin{center}
\includegraphics[width=6cm]{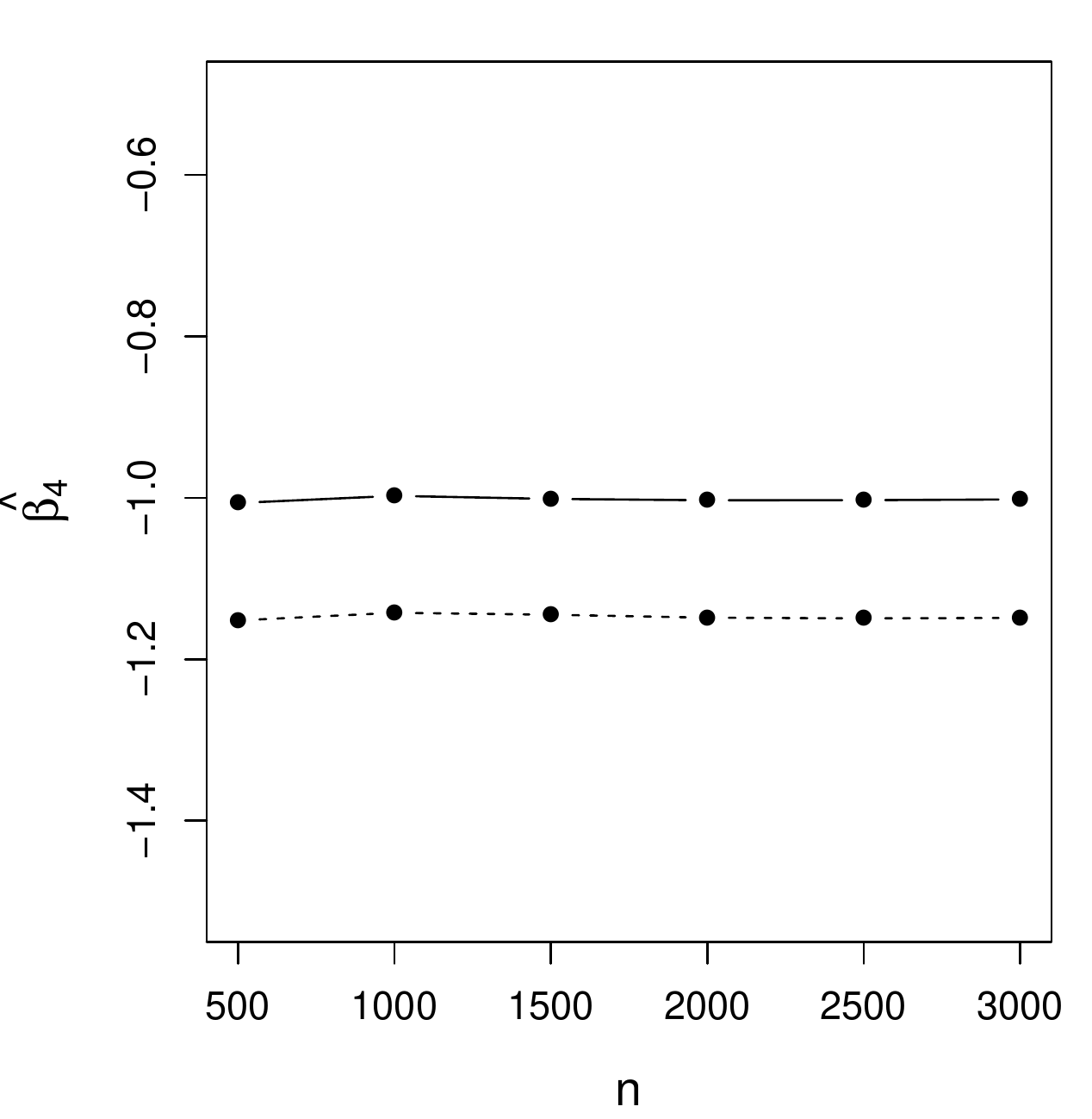}
\includegraphics[width=6cm]{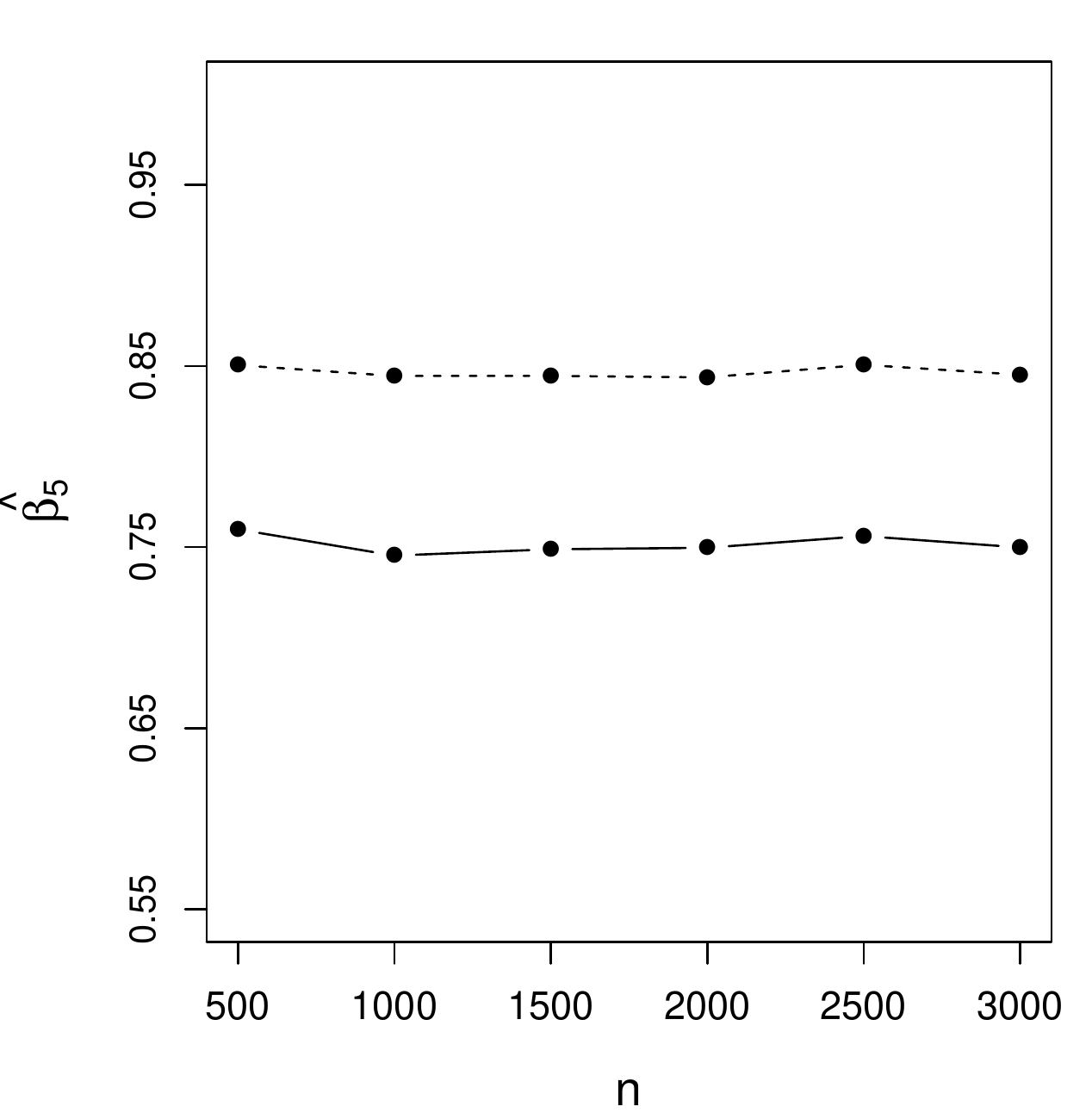}\\
\includegraphics[width=6cm]{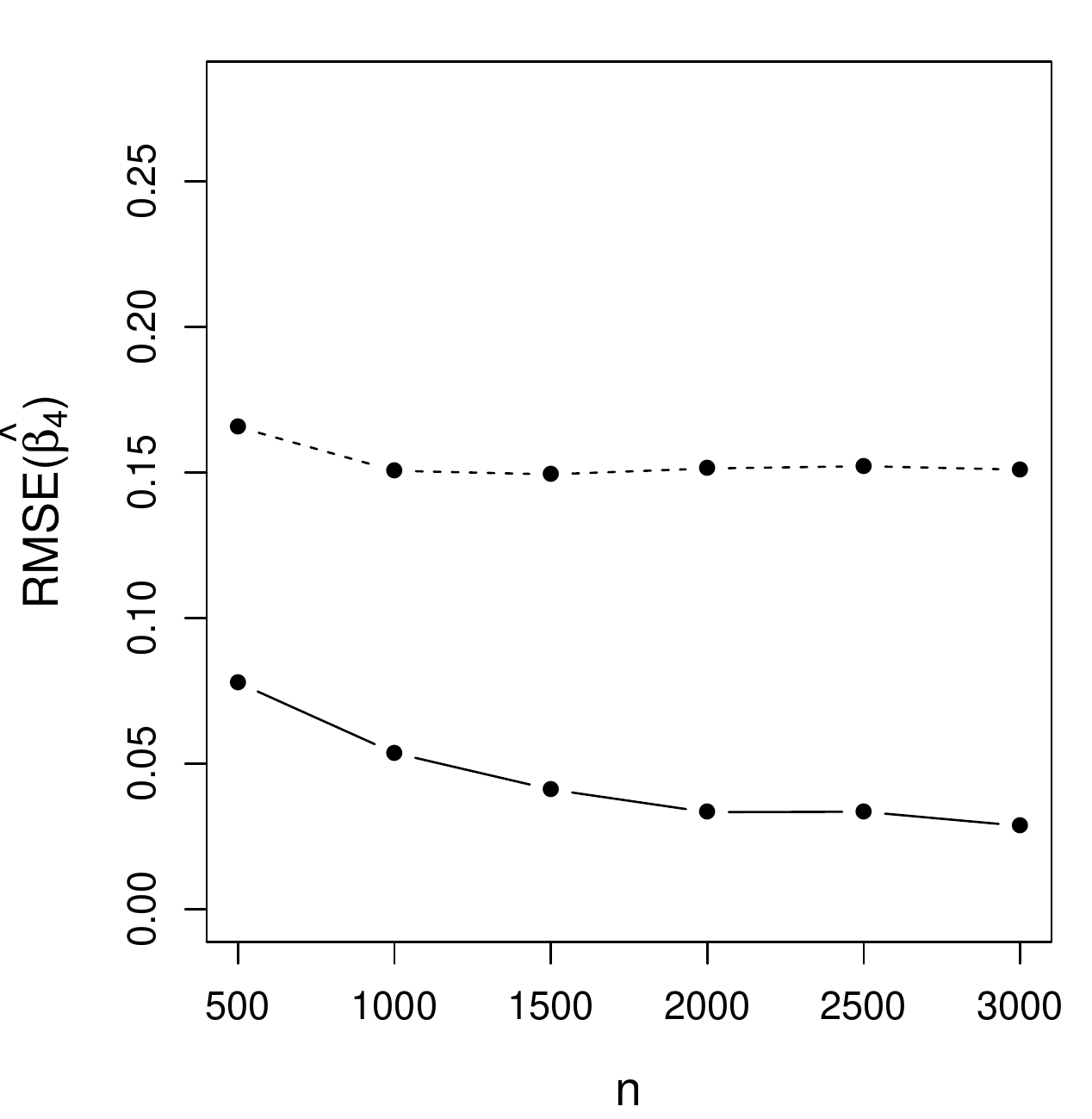}
\includegraphics[width=6cm]{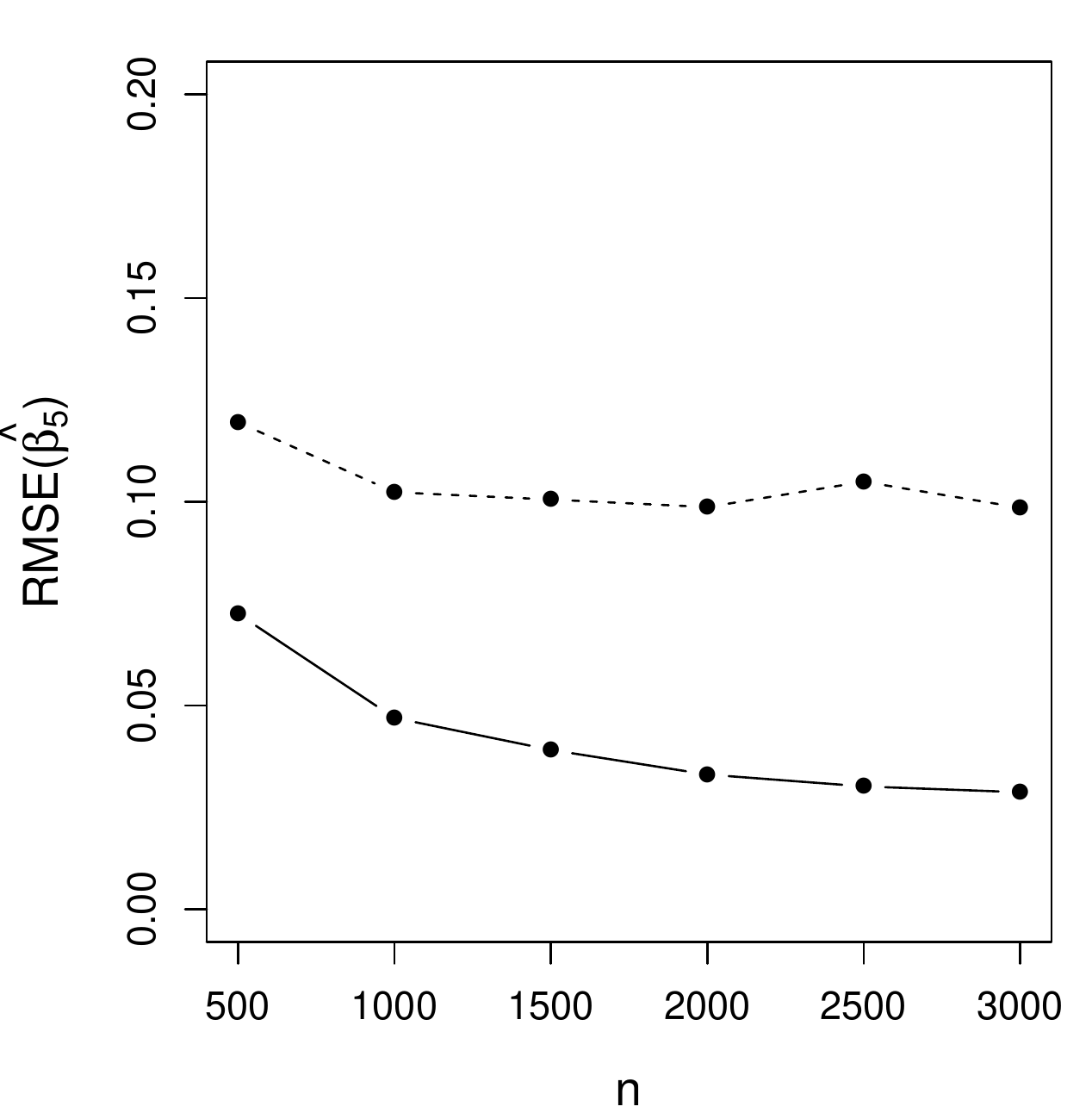}
\end{center}
\caption{Empirical expected values of the GASSM and GAM estimators for $\beta_4=-1$ and $\beta_5=0.75$ (top plots) and square root of their mean squared errors (bottom plots). GASSM denotes the generalized additive sample selection model while GAM denotes the classic univariate generalized additive model. Solid lines refer to GASSM results whereas dotted lines to those produced using GAM.}
\label{fig.s1}
\end{figure}

\begin{figure}
\begin{center}
\includegraphics[width=6cm]{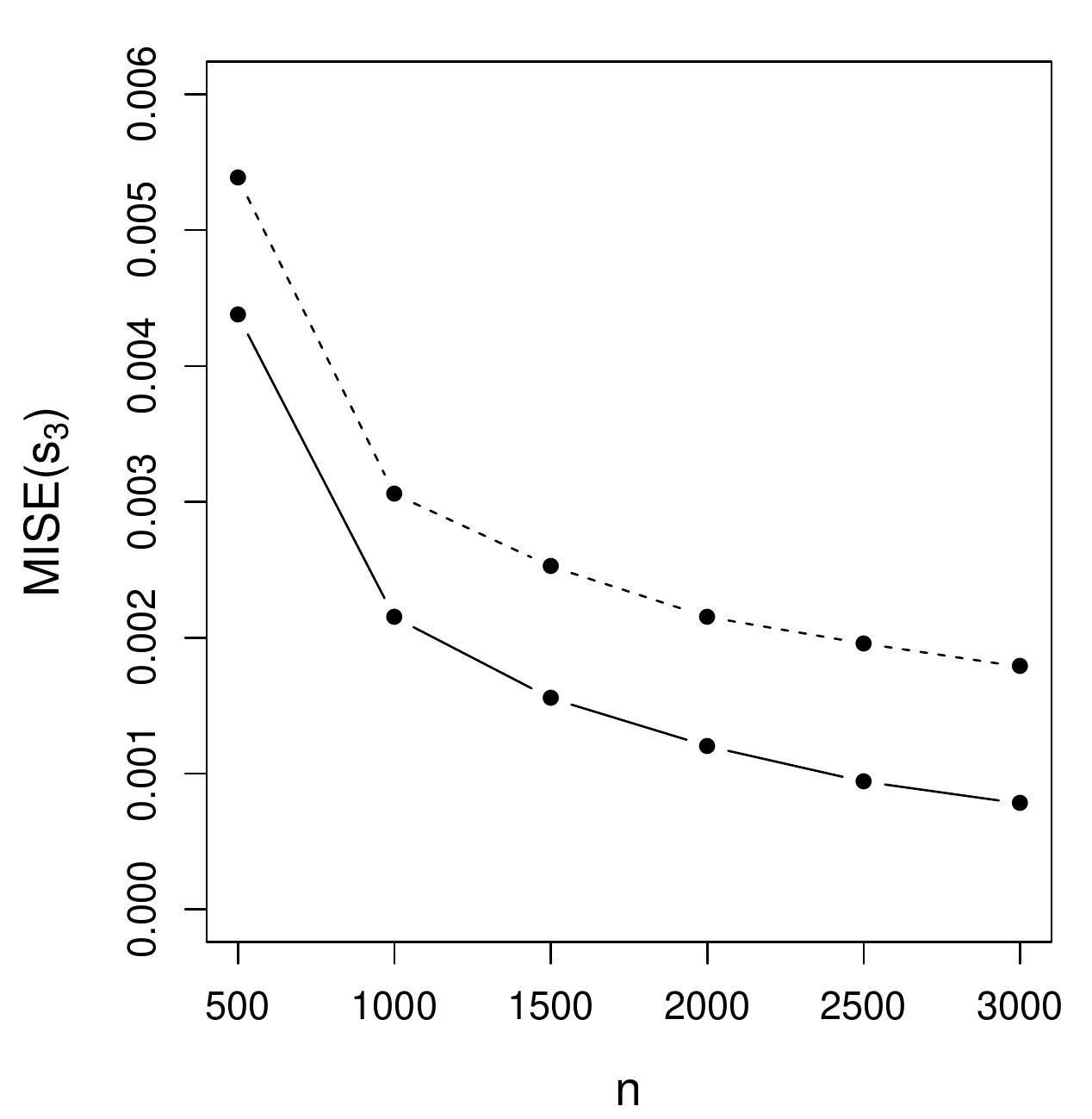}
\includegraphics[width=6cm]{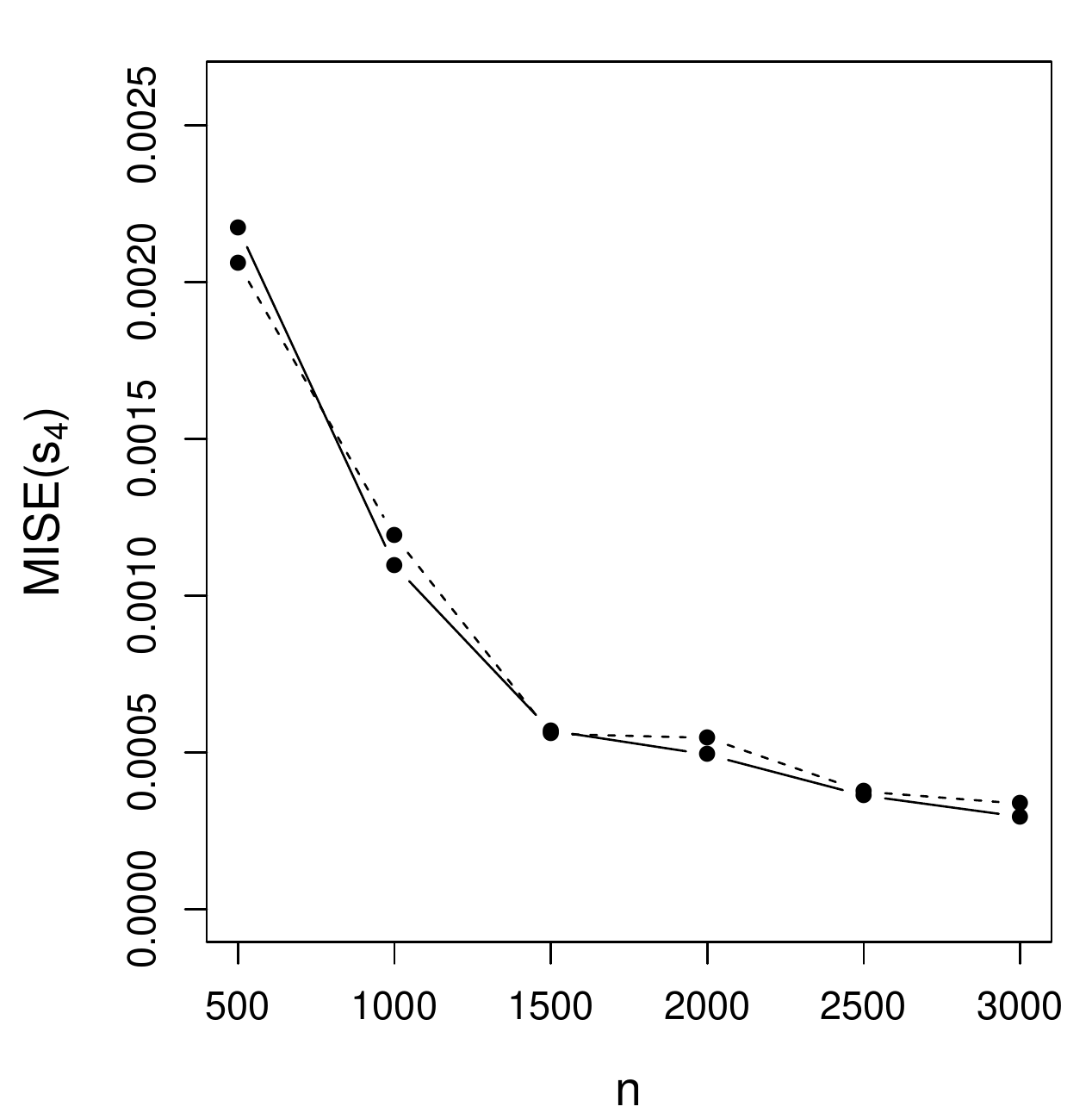}
\end{center}
\caption{Empirical mean integrated squared errors of the GASSM and GAM estimators for $s_3$ and $s_4$.}
\label{fig.s3}
\end{figure}

The empirical expected values of the GASSM and GAM estimators for $\beta_4$ and $\beta_5$ and the square root of their empirical mean squared errors are plotted as a function of sample size and shown in Figure \ref{fig.s1}. Moreover, the mean integrated squared errors (MISE) of the estimators for $s_3$ and $s_4$ are shown in Figure \ref{fig.s3}. In all plots, solid lines correspond to the GASSM results and dotted lines to those obtained using GAM. The two top plots in Figure \ref{fig.s1} show that GAM estimators are biased, while the GASSM estimators have empirical mean values practically identical to the true values. It is worth noting that the GAM estimators have smaller standard deviations than those of GASSM (see Table \ref{Table.s1}). This is because the latter acknowledge the uncertainty due the selection mechanism. Importantly, the mean squared errors of the GASSM estimators remain uniformly smaller than those for GAM as shown in the two bottom plots of Figure \ref{fig.s1}.

Empirical consistency of the estimators for $s_3$ and $s_4$ can be studied by looking at the results in Figure \ref{fig.s3}. Interestingly, while the GAM estimator for $s_3$ has a MISE which is considerably larger than that of the GASSM, both estimators for $s_4$ appear to be equally good; a possible interpretation is that function $s_4$ is linear and the additional parameters associated with its basis functions available for estimation may compensate for selection bias.

\begin{table}
\begin{center}
\begin{tabular}{ccccccc}
\hline\hline
$n$   & $500$    & $1000$   & $1500$                  & $2000$   & $2500$   & $3000$ \\\hline
      &          &          & ${\rm SD}(\hat\beta_4)$ &          &          &     \\
GASSM &	$0.0778$ & $0.0538$ & $0.0413$                & $0.0334$ & $0.0335$ & $0.0287$\\
GAM   & $0.0663$ & $0.0492$ & $0.0370$                & $0.0298$ & $0.0293$ & $0.0263$\\\hline
      &          &          & ${\rm SD}(\hat\beta_5)$ &          &          &     \\
GASSM &	$0.0720$ & $0.0469$ & $0.0392$                & $0.0332$ & $0.0296$ & $0.0288$	\\
GAM   & $0.0650$ & $0.0395$ & $0.0344$                & $0.0311$ & $0.0298$ & $0.0258$  \\\hline\hline
\end{tabular}
\end{center}
\caption{Empirical standard deviations of the GASSM and GAM estimators of $\beta_4=-1$ and $\beta_5=0.75$.}
\label{Table.s1}
\end{table}

\subsection{Comparison to logged model}

In real world applications, the analysis of data with positive outcomes having a highly skewed distribution is often performed using log-transformed outcomes. Then a model assuming a normal or $t$ distribution is fitted (see, e.g., the example in \cite{Marchenko-and-Genton-2012}). This section shows the potential inaccuracy of this practice.

For the latent selection variable, it was assumed that
$$
Y_{1i}^*=\alpha_0 + \alpha_1 x_1 + \alpha_2 x_2 + \alpha_3 x_3 + \varepsilon_i,
$$
with $\alpha_0=0.58$, $\alpha_1=2.5$, $\alpha_2=-1$, $\alpha_3=0.8$ and $\varepsilon_i\sim N(0,1)$, whereas the outcome variable $Y_{2i}^*$ was assumed to pertain to gamma distribution with shape parameter $k=2$ and expected value $\mu_i=\E( Y_{2i}^*)$ such that 
$$
 \log( \mu_i ) = \beta_0 + \beta_1 x_1 + \beta_2 x_2,
$$
with $\beta_0=-0.68$, $\beta_1=-1.5$ and $\beta_2=0.5$. Three different patterns of dependence between $Y_1^*$ and $Y_2^*$ specified by the normal, Frank and Clayton copulae were considered. For each copula, three values of association parameter $\theta$, corresponding to the values of Kendall's $\tau$ equal to $0.1$, $0.5$ and $0.7$, were used.

To generate covariates $x_1$, $x_2$ and $x_3$, random numbers $(z_1,z_2,z_3)$ pertaining to a trivariate normal distribution were used. That is, we first generated  
$$
(z_1,z_2,z_3)\sim N\left( 
\left[ \begin{array}{c}
0\\
0\\
\end{array} \right]
, \,
\left[ \begin{array}{ccc}
1 & 0.5 & 0.5 \\
0.5 & 1 & 0.5 \\
0.5 & 0.5 & 1\\
\end{array} \right]
 \right)
$$
and then $x_1$, $x_2$ and $x_3$ were obtained as $x_1=\mathbbm{1}(\Phi^{-1}(z_1)>0.5)$, $x_2=\Phi^{-1}(z_2)$ and $x_3=\Phi^{-1}(z_3)$. Thus, predictor $x_1$ was a binary variable, whereas $x_1$ and $x_2$ were uniformly distributed on $(0,1)$, with correlation approximately equal to 0.5.

For each generated data set, two models were fitted: one assuming a normal distribution for the logarithm of the outcome and the other assuming a gamma distribution for the outcome. The sample size was $n=1000$ and the number of Monte Carlo repetitions $300$. For both fitted models, the relative bias and root of mean square error of estimators $\hat\beta_0$, $\hat\beta_1$, $\hat\beta_2$ and Kendall's $\hat\tau$ (related to $\hat\theta$) were calculated and are reported in Table \ref{lognorm}, which also shows test errors of the fitted models.

\begin{table}
{\footnotesize
\begin{tabular}{cccccccccccc}
\hline \hline
 &  & \multicolumn{2}{c}{$\hat\beta_0$} &  \multicolumn{2}{c}{$\hat\beta_1$} &  \multicolumn{2}{c}{$\hat\beta_2$}  & \multicolumn{2}{c}{$\hat\tau$} & Test \\[2pt]
 &  & Bias (\%) &  RMSE & Bias (\%) &  RMSE & Bias (\%) &  RMSE & Bias (\%) &  RMSE &   error \\[5pt] 
 \hline
  \multicolumn{11}{c}{ Normal Copula} \\[5pt]
 \hline
\multirow{2}{*}{\rotatebox{90}{$\tau=0.1$}} 
 & G 	 & 	 -5.7 	 & 	 0.129 	 & 	 3.8 	 & 	 0.15 	 & 	 5.3 	 & 	 0.105 	 & 	 -83.5 	 & 	 0.249 	 & 	 0.329 \\[1pt]
 & L 	 & 	 23.7 	 & 	 0.253 	 & 	 9.7 	 & 	 0.272 	 & 	 11.9 	 & 	 0.141 	 & 	 -227.4 	 & 	 0.426 	 & 	 0.339 \\[10pt]

 \multirow{2}{*}{\rotatebox{90}{$\tau=0.5$}}
 & G 	 & 	 -0.4 	 & 	 0.069 	 & 	 0.7 	 & 	 0.077 	 & 	 1.9 	 & 	 0.086 	 & 	 -1 	 & 	 0.13 	 & 	 0.321 \\[1pt]
 & L 	 & 	 30.9 	 & 	 0.236 	 & 	 5.3 	 & 	 0.153 	 & 	 5.9 	 & 	 0.111 	 & 	 -19.3 	 & 	 0.243 	 & 	 0.334 \\[10pt]
 
 \multirow{2}{*}{\rotatebox{90}{$\tau=0.7$}}
 & G 	 & 	 -0.1 	 & 	 0.06 	 & 	 0.5 	 & 	 0.066 	 & 	 1.9 	 & 	 0.084 	 & 	 0.3 	 & 	 0.104 	 & 	 0.321 \\[1pt]
 & L 	 & 	 32.3 	 & 	 0.229 	 & 	 4.2 	 & 	 0.096 	 & 	 4.1 	 & 	 0.098 	 & 	 -7.1 	 & 	 0.121 	 & 	 0.332 \\[5pt]
 
  \hline
  \multicolumn{11}{c}{ Frank Copula} \\[5pt]
 \hline
\multirow{2}{*}{\rotatebox{90}{$\tau=0.1$}} 
 & G 	 & 	 -6.1 	 & 	 0.13 	 & 	 3.2 	 & 	 0.148 	 & 	 0.5 	 & 	 0.094 	 & 	 -34 	 & 	 0.245 	 & 	 0.328 \\[1pt]
 & L 	 & 	 39.7 	 & 	 0.32 	 & 	 -0.4 	 & 	 0.204 	 & 	 -3.5 	 & 	 0.12 	 & 	 18.9 	 & 	 0.318 	 & 	 0.347 \\[10pt]
 \multirow{2}{*}{\rotatebox{90}{$\tau=0.5$}}
  & G 	 & 	 -2.7 	 & 	 0.085 	 & 	 1.4 	 & 	 0.095 	 & 	 -0.3 	 & 	 0.084 	 & 	 -6.8 	 & 	 0.18 	 & 	 0.324 \\[1pt]
 & L 	 & 	 33 	 & 	 0.249 	 & 	 3.2 	 & 	 0.132 	 & 	 -0.6 	 & 	 0.102 	 & 	 -5.1 	 & 	 0.205 	 & 	 0.338 \\[10pt]

 \multirow{2}{*}{\rotatebox{90}{$\tau=0.7$}}
  & G 	 & 	 -1.6 	 & 	 0.07 	 & 	 0.8 	 & 	 0.078 	 & 	 -0.3 	 & 	 0.084 	 & 	 -2.5 	 & 	 0.149 	 & 	 0.324 \\[1pt]
 & L 	 & 	 30.8 	 & 	 0.225 	 & 	 4.2 	 & 	 0.112 	 & 	 0.1 	 & 	 0.098 	 & 	 -5.9 	 & 	 0.156 	 & 	 0.335 \\[5pt]

   \hline
  \multicolumn{11}{c}{ Clayton Copula} \\[5pt]
 \hline
\multirow{2}{*}{\rotatebox{90}{$\tau=0.1$}} 
 & G 	 & 	 1.3 	 & 	 0.058 	 & 	 -0.4 	 & 	 0.064 	 & 	 1.9 	 & 	 0.094 	 & 	 5.5 	 & 	 0.071 	 & 	 0.322 \\[1pt]
 & L 	 & 	 32.3 	 & 	 0.227 	 & 	 3.9 	 & 	 0.086 	 & 	 4.9 	 & 	 0.107 	 & 	 -76.5 	 & 	 0.085 	 & 	 0.332 \\[10pt]
 \multirow{2}{*}{\rotatebox{90}{$\tau=0.5$}}
 & G 	 & 	 0.6 	 & 	 0.047 	 & 	 0 	 & 	 0.054 	 & 	 1.4 	 & 	 0.086 	 & 	 1.7 	 & 	 0.06 	 & 	 0.321 \\[1pt]
 & L 	 & 	 32.4 	 & 	 0.229 	 & 	 4.1 	 & 	 0.091 	 & 	 4.7 	 & 	 0.101 	 & 	 -7.9 	 & 	 0.093 	 & 	 0.332 \\[10pt]
 \multirow{2}{*}{\rotatebox{90}{$\tau=0.7$}}
 & G 	 & 	 0.4 	 & 	 0.045 	 & 	 0.1 	 & 	 0.05 	 & 	 1.1 	 & 	 0.085 	 & 	 1.4 	 & 	 0.052 	 & 	 0.32 \\[1pt]
 & L 	 & 	 33.5 	 & 	 0.233 	 & 	 3.6 	 & 	 0.077 	 & 	 4.1 	 & 	 0.097 	 & 	 -1.8 	 & 	 0.053 	 & 	 0.332 \\[5pt]
 \hline \hline
\end{tabular}
}
\caption{Comparison of performance between the gamma sample selection model (G) and the normal sample selection model in which the logarithm of the outcome variable is used (L).}
\label{lognorm}
\end{table}

In most cases considered the gamma sample selection model outperforms the normal one (which employes a log-transformed outcome) in terms of bias, mean squared error and test error. This not only shows that the proposed model is flexible enough to accommodate non-Gaussian distributions, but also that using a transformed outcome can lead to unreliable empirical results.

\section{Discussion}\label{discussion}

We have introduced an extension of the generalized additive model which accounts for non-random sample selection. The proposed approach is flexible in that it allows for different distributions of the outcome variable, several dependence structures between the outcome and selection equations, and non-parametric effects on the responses. Parameter estimation with integrated automatic multiple smoothing parameter selection is achieved within a penalized likelihood and simultaneous equation framework. We have established the asymptotic theory for the proposed penalized spline estimators, and illustrated the empirical effectiveness of the approach through a simulation study. A few points are noteworthy.

\begin{itemize}
\item The generalized sample selection model has been formulated using penalized B-splines. This allows for simple handling of the model's theoretical properties. However, in practice different smoothers can be used, for example truncated polynomials (which yield an equivalent approach as detailed in \cite{Krivobokova-et-al-2009}) or thin plate regression splines \citep{Wood} .

\item The estimation procedure discussed in Section \ref{Estimation} has been implemented in the freely available {\tt R} package {\tt SemiParSampleSel}. Currently, the outcome can be modeled using the normal, gamma and a number of discrete distributions. The copulae available are: normal, Clayton, Joe, Frank, Gumbel, AMH, FGM and their rotated versions. Given the modular structure of the estimation algorithm, other copulae and marginal distributions can be incorporated in {\tt SemiParSampleSel} with little programming work. 

\item For simplicity of treatment and notation, the generalized additive sample selection model has been defined using as few parameters as possible. However, many model structures are allowed within the proposed framework. For example, the association parameter $\theta$ can be made dependent on predictors and hence enter the likelihood function as a transformed linear predictor instead of a scalar. Similarly, even though the scale parameter $\phi$ has been set to $1$ for simplicity of derivations, additional parameters related to the specific distribution employed can be estimated. In both cases, all the theoretical derivations presented in the paper still hold.

\item Assumption {\bf A4} in Section \ref{Asymptotics} allows the sequence smoothing parameters $\hat{\boldsymbol\lambda}_n$ to grow as the sample size increases. This condition is rather week as, in fact, the sequence $\hat{\boldsymbol\lambda}_n$ based on the mean squared error criterion described in Section A.2 is bounded in probability (cf e.g., Kauermann, 2005). Thus the theoretical properties of the penalized estimator derived in Section \ref{Asymptotics} hold even if the smoothing parameters are estimated, and not deterministic.

\end{itemize}

An interesting direction of future research will be to compare the small sample performances of the proposed estimator and some of the non-parametric, semiparametric and Bayesian methods mentioned in the introduction. Moreover, for many copulae a specific value of the association parameter $\theta$ yields a product distribution which indicates lack of non-random sample selection. Thus, the important issue of testing hypotheses regarding parameter $\theta$ will be addressed in future work.
\\

\noindent
{\bf Acknowledgements}\\
\noindent
This research was supported by the Engineering and Physical Sciences Research Council, UK (Grant EP\/J006742\/1).

\begin{appendix}
\noindent

\section{Algorithmic details}
\label{Estimation_details}

\subsection{Trust region algorithm}
\label{trust}

Recall that $\bm\delta=(\boldsymbol\alpha,\boldsymbol\beta,\theta)$ and define the penalized gradient and Hessian at iteration $a$ as $G^{[a]}_p=G^{[a]}-\textbf{S}_{\hat{\bm\lambda}}\bm\delta^{[a]}$ and $H^{[a]}_p=H^{[a]}-\textbf{S}_{\hat{\bm\lambda}}$. Each iteration of the trust region algorithm solves the problem
\begin{equation*}
\begin{aligned}
\underset{\textbf{p}}{\operatorname{min}} \ \ \breve{\ell_p}(\bm\delta^{[a]})&\defeq-\left\{\ell_p(\bm\delta^{[a]})+\textbf{p}\ts{G}^{[a]}_p+\frac{1}{2}\textbf{p}\ts{H}^{[a]}_p\textbf{p}\right\} \ \ \text{such that} \ \ \|\textbf{p}\|\leq r^{[a]},\\
\bm\delta^{[a+1]}&=\underset{\textbf{p}}{\operatorname{arg \ min}} \ \ \breve{\ell_p}(\bm\delta^{[a]})+\bm\delta^{[a]},
\end{aligned}
\end{equation*}
where $\|\cdot\|$ denotes the Euclidean norm, and $r^{[a]}$ is the radius of the trust region. At each iteration of the algorithm, $\breve{\ell_p}(\bm\delta^{[a]})$ is minimized subject to the constraint that the solution falls within a trust region with radius $r^{[a]}$. The proposed solution is then accepted or rejected and the trust region expanded or shrunken based on the ratio between the improvement in the objective function when going from $\bm\delta^{[a]}$ to $\bm\delta^{[a+1]}$ and that predicted by the quadratic approximation. See \cite{Geyer} for the exact details (e.g., numerical stability and termination criteria) of the implementation used here. It is important to stress that near the solution the trust region method typically behaves as a classic unconstrained algorithm \citep{Geyer,Nocedal}. Starting values for the coefficients in $\bm\alpha$ and $\bm\beta$ are obtained by fitting the selection and outcome equations separately. The initial parameter of $\theta$ is set to zero as there is not typically good {\it a priori} information about the direction and strength of the association between the selection and outcome equations, conditional on covariates.

\subsection{Multiple smoothing parameter estimation}\label{smoothness} 
\label{estimating_lambda}


Let us use the fact that near the solution the trust region algorithm usually behaves as a classic Newton or Fisher Scoring method, and assume that $\bm\delta^{[a+1]}$ is a new updated guess for the parameter vector which maximizes $\ell_p$. If $\bm\delta^{[a+1]}$ is to be `correct', then the penalized gradient evaluated at those parameter values would be $\textbf{0}$, i.e. $G_p^{[a+1]}=\textbf{0}$. Applying a first order Taylor expansion to $G^{[a+1]}_p$ about $\bm\delta^{[a]}$ yields
$\textbf{0}=G_p^{[a+1]} \approx G_p^{[a]} + \left( \bm\delta^{[a+1]}-\bm\delta^{[a]} \right)H^{[a]}_p$, from which we find the solution at iteration $a+1$. After some manipulation, this can be expressed as 
$$
\bm\delta^{[a+1]}= \left(\In^{[a]} + \Sl\right)^{-1} \sqrt{\In^{[a]}} \mathbf{z}^{[a]},
$$
where $\In^{[a]}$ is $-H^{[a]}$ (or, alternatively, $-\E\left(H^{[a]}\right)$), and $\mathbf{z}^{[a]}=\sqrt{\In^{[a]}}\bm\delta^{[a]} + \eb^{[a]}$, with $\eb^{[a]}=\sqrt{\In^{[a]}}^{-1}G^{[a]}$. From standard likelihood theory, $\eb \sim \mathcal{N}\left(\textbf{0},\I \right)$ and $\mathbf{z} \sim \mathcal{N} \left(\muz, \I  \right)$, where $\I$ is an identity matrix, $\muz=\sqrt{\In}\bm\delta^0$, and $\bm\delta^0$ is the true parameter vector. The predicted value vector for $\mathbf{z}$ is $\hat{\bm\mu}_{\textbf{z}}=\sqrt{\In}\hat{\bm\delta}=\textbf{A}_{\hat{\bm\lambda}}\mathbf{z}$, where $\textbf{A}_{\hat{\bm\lambda}}=\sqrt{\In}\left(\In + \Sl\right)^{-1} \sqrt{\In}$. Since our goal is to select the smoothing parameters in as parsimonious manner as possible so that the smooth terms' complexity which is not supported by the data is suppressed, $\bm\lambda$ is estimated so that $\hat{\bm\mu}_{\textbf{z}}$ is as close as possible to $\muz$. This can be achieved using 
\begin{equation}
\begin{split}
\E \left( \| \muz - \hat{\bm\mu}_{\textbf{z}} \|^2 \right) & = \E \left(\| \mathbf{z} - \textbf{A}_{\bm\lambda}\mathbf{z} - \eb    \|^2  \right)\\
& = \E\left( \| \mathbf{z} - \Al\mathbf{z} \|^2\right)+ \E\left(- \eb\ts\eb -2\eb\ts\muz + 2\eb\ts\Al\muz+2\eb\ts\Al\eb\right)\\
 & = \E\left( \| \mathbf{z} - \Al\mathbf{z} \|^2\right) - \check{n} + 2\text{tr}(\Al),
  \end{split}
  \nonumber
\end{equation}
where $\check{n}=3n$ and $\text{tr}(\textbf{A}_{\bm\lambda})$ is the number of effective degrees of freedom of the penalized model. Hence, the smoothing parameter vector is estimated by minimizing an estimate of the expectation above, that is
\beq\label{AICc}
\mathcal{V}(\bm\lambda)= \| \mathbf{z} - \Al\mathbf{z} \|^2 - \check{n} + 2\text{tr}(\Al),
\eeq 
which is equivalent to the expression of the Un-Biased Risk Estimator given in \citet[][Chapter 4]{Wood}. This is also equivalent to the Akaike information criterion after dropping the irrelevant constant; the first term on the right hand side of (\ref{AICc}) is a quadratic approximation to $-2\ell(\hat{\bm\delta})$ to within an additive constant. In practice, given $\bm\delta^{[a+1]}$, the problem becomes 
\beq\label{PIRLS1}
{\bm\lambda}^{[a+1]}=\underset{\bm\lambda}{\operatorname{arg \ min}} \ \ \mathcal{V}(\bm\lambda)\defeq\| \mathbf{z}^{[a+1]} - \Al^{[a+1]}\mathbf{z}^{[a+1]} \|^2 - \check{n} + 2\text{tr}(\Al^{[a+1]}),
\eeq
which is solved using the automatic stable and efficient computational routine by \cite{Wood04}.

\section{Proofs of Lemmas}
\label{proofs}

{\bf Proof of Lemma 2.}
We have 
$$
F_{n}(\bdelta^0) = 
	  \left[ \begin{array}{ccc} \left({\bf X}^{(1)}\right)^T & {\bf 0} & {\bf 0} \\ {\bf 0} &  \left({\bf X}^{(2)}\right)^T & {\bf 0} \\
		                          {\bf 0} & {\bf 0} & {\bf 1}^T   \end{array} \right]
	  \left[ \begin{array}{ccc} \E_{\bdelta^0}(W_1) & \E_{\bdelta^0}(W_3) & \E_{\bdelta^0}(W_4) \\
		                          \E_{\bdelta^0}(W_3) &  \E_{\bdelta^0}(W_2) & \E_{\bdelta^0}(W_5) \\
															\E_{\bdelta^0}(W_4) &  \E_{\bdelta^0}(W_5) & \E_{\bdelta^0}(W_6) \end{array} \right]
	  \left[ \begin{array}{ccc} {\bf X}^{(1)} & {\bf 0} & {\bf 0} \\ {\bf 0} &  {\bf X}^{(2)} & {\bf 0} \\
  	                          {\bf 0} & {\bf 0} & {\bf 1} \end{array} \right].
$$
From Lemma 1 of \cite{Yoshida-Naito-2012} we obtain that elements of matrices $\left(X_j^{(1)}\right)^T X_j^{(1)}$ are of order $O\left(\frac{n}{K_n}\right)$ for $j=1,\ldots,D_1$, and elements of matrices $\left(X_j^{(1)}\right)^T X_l^{(1)}$ are of order $O\left(\frac{n}{K_n^2}\right)$ for $j\neq l$, $j,l=1,\ldots,D_1$. The same boundaries hold for the matrices $X_j^{(2)}$, $j=1,\ldots,D_2$.
\\
Thus elements of matrices $\left({\bf X}^{(1)}\right)^T {\bf X}^{(1)}$ and $\left({\bf X}^{(2)}\right)^T {\bf X}^{(2)}$ are of order $O\left(\frac{n}{K_n}\right)$. In a straightforward way, the result also extends to the matrix $\left({\bf X}^{(1)}\right)^T {\bf X}^{(2)}$.
\\
Now we consider the order of the diagonal elements $w_i^{(1)}$, \ldots, $w_i^{(6)}$. It holds
\begin{equation}
\left\lvert \frac{\partial F_{2i}}{\partial\eta_{2i}}(y_2) \right\rvert  = \left\lvert  \int_{-\infty}^{y_2} (v-b'(\bbeta^0)) f_{2i}(v) dv \right\rvert  
 \leq \E_{\bdelta^0}|Y_{2i}| + |b'(\bbeta^0)| \leq 2\E_{\bdelta^0}|Y_{2i}|.
\label{pr1}
\end{equation}
and
\begin{equation}
  \frac{\partial^2 F_{2i}}{\partial\eta_{2i}^2}(y_2) = \int_{-\infty}^{y_2} \left( 1-b''(\eta_{2i}) + (v-b'(\eta_{2i}))^2 \right) f_{2i}(v) dv.
\label{pr2}
\end{equation}
Thus $|\frac{\partial^2 F_{2i}}{\partial\eta_{2i}^2}(y_2)| \leq 1+2 {\rm Var}_{\bdelta^0}(Y_{2i})$. This combined with (\ref{pr1}) and assumptions (A1) and (A2) yields $\E(w_i^{(j)}) = O(1)$ for $j=1,\ldots,6$.

\vspace{2mm}\noindent
Moreover, by the properties of B-spline basis the $(i,l)$th components of $\left(X_j^{(1)}\right)^T X_k^{(1)}$, for $j,k=1,\ldots,D_1$, and $\left(X_j^{(2)}\right)^T X_k^{(2)}$, for $j,k=1,\ldots,D_2$, equal 0 if $|i-l|>p$. Hence the matrices $\left({\bf X}^{(1)}\right)^T W_1{\bf X}^{(1)}$, $\left({\bf X}^{(2)}\right)^T W_2{\bf X}^{(2)}$ are band matrices and the assertion follows.
\begin{flushright}
\vspace{-5mm}
$\square$
\end{flushright}

\noindent
{\bf Proof of Lemma 3.} 
We use induction w.r.t. the number of variables. Let $M_n=D_1+D_2+1$ and matrix $U_{M_n}$ be defined as
$$
  U_{M_n} = F_{n,p}(\bdelta^0) =
\alpha	 \left[ 
	\begin{array}{cc}
	 U_{M_n-1} & R^T\\
	 R & \Lambda
	\end{array}
	\right].
$$
The result of \cite{Horn-Johnson-1985} yields
$$
  U_{M_n}^{-1} = 
		\left[ 
	\begin{array}{cc}
	 U_{M_n-1}^{-1}+U_{M_n-1}^{-1}R^TV^{-1}RU_{M_n-1}^{-1} & -U_{M_n-1}^{-1}R^TV^{-1}\\
	 -V^{-1}RU_{M_n-1}^{-1}  & U_{M_n-1}^{-1}
	\end{array}
	\right],
$$
where $V^{-1}=\Lambda-RU_{M_n-1}^{-1}R^T$.
Then assertion can be proven similarly to \cite{Krivobokova-et-al-2009} by using the fact that matrices $\left({\bf X}^{(1)}\right)^T W_1 {\bf X}^{(1)}$, $\left({\bf X}^{(2)}\right)^T W_1 {\bf X}^{(2)}$ and $\left({\bf X}^{(1)}\right)^T W_1 {\bf X}^{(2)}$ are band matrices and the properties of the inverse of band matrices listed in \cite{Demko-1977}.
\begin{flushright}
\vspace{-5mm}
$\square$
\end{flushright}

\noindent
{\bf Proof of Lemma 4} (sketch). 
$$
H_n(\bdelta^0) - F_n(\bdelta^0) = {\bf X}^T (W - \E(W)) {\bf X}. 
$$
It holds $w_i-\E(w_i)=O_P(n^{-1/2})$ as every $w_i$ is a sum of independent and bounded random variables. Moreover,
$$
  \frac{1}{n} \sum_{i=1}^n \left(B_{-p+j}(x_i) B_{-p+l}(x_i)\right)^2 = O(K_n^{-1}).
$$
Hence $\sum_{i=1}^n{\rm Var}(w_i B_{-p+j}(x_i)B_{-p+l}(x_i)) = O(n/K_n)$ which yields the assertion.
\begin{flushright}
\vspace{-5mm}
$\square$
\end{flushright}

\end{appendix}

\bibliographystyle{apalike2}
\bibliography{ref}

\end{document}